\documentclass[11pt]{article}
\usepackage{eurosym}
\usepackage{amsmath, amssymb}
\usepackage{graphicx}
\usepackage{color}

\DeclareMathOperator{\dom}{dom}

\DeclareMathOperator{\gph}{gph}
\DeclareMathOperator{\clm}{clm}
\DeclareMathOperator{\Lipusc}{Lipusc}
\DeclareMathOperator{\Hof}{Hof}

\DeclareMathOperator{\extr}{extr}
\DeclareMathOperator{\spann}{span}

\DeclareMathOperator{\cone}{cone}

\newtheorem{theo}{Theorem}
\newtheorem{lem}{Lemma}
\newtheorem{prop}{Proposition}
\newtheorem{cor}{Corollary}
\newtheorem{rem}{Remark}
\newtheorem{exa}{Example}
\newtheorem{defn}{Definition}

\newenvironment{dem}[1][Proof]{\noindent \textbf{#1.} }{\ \rule{0.5em}{0.5em}}

\input{tcilatex}
\begin{document}

\title{Hoffman constant of the argmin mapping in linear optimization \thanks{%
This research has been partially supported by Grant PID2022-136399NB-C22
from MICINN, Spain, and ERDF, "A way to make Europe", European Union, and
Grant PROMETEO/2021/063 from Generalitat Valenciana, Spain. This research
was also partially supported by the CIPROM/2024/34 grant, funded by the
Conselleria de Educaci\'{o}n, Cultura, Universidades y Empleo, Generalitat
Valenciana.}}
\author{J. Camacho\thanks{%
Center of Operations Research, Miguel Hern\'{a}ndez University of Elche,
03202 Elche (Alicante), Spain (j.camacho@umh.es, canovas@umh.es,
parra@umh.es).} \and M.J. C\'{a}novas\footnotemark[2] \and H. Gfrerer\thanks{%
Johann Radon Institute for Computational and Applied Mathematics (RICAM),
A-4040 Linz, Austria and Institute of Information Theory and Automation,
Czech Academy of Sciences, 18208 Prague, Czech Republic;
(helmut.gfrerer@ricam.oeaw.ac.at).} \and J. Parra\footnotemark[2]}
\date{}
\maketitle

\begin{abstract}
The main goal of this paper is to provide a point-based expression for the
Hoffman constant of the argmin mapping in linear optimization, understood as
the sharp Lipschitz constant restricted to its domain. The work is mainly
developed in the parametric context of right-hand side perturbations of the
constraint system. To the authors' knowledge, this is the first exact
formula for this constant, although we can find in the literature different
upper estimates. The paper tackles this objective from a broader
perspective, which introduces new tools of their own interest, such as the
concept of well-connected piecewise convex mapping. We isolate the nice
behavior of such mappings to derive a crucial equality between the Hoffman
constant (which is a global stability measure) and the supremum of calmness
moduli (of local nature). The paper also includes some specifics about
directional stability of optimal solutions and finishes with some
conclusions and notes about further research.


\bigskip

\end{abstract}

\section{Introduction}

The primary purpose of this paper is to derive the \emph{Hoffman constant}
for the optimal solutions of linear optimization problems given in the form%
\begin{equation}
\begin{array}{ll}
\text{minimize } & c^{\prime }x \\ 
\text{subject to } & Ax\leq b,%
\end{array}
\label{eq_LPproblem}
\end{equation}%
where $x\in \mathbb{R}^{n}$ is the decision variable, $A\in \mathbb{R}%
^{m\times n}$ is a given matrix, and $c\in \mathbb{R}^{n}$ and $b\in \mathbb{%
R}^{m}$ are regarded as parameters. Elements in $\mathbb{R}^{n}$ are
considered as column vectors and the prime stands for transposition; i.e., $%
y^{\prime }$ denotes the transpose of $y\in \mathbb{R}^{n}.$

Throughout the paper, matrix $A$ is fixed and we are focused on the \emph{%
optimal set (argmin) }mapping, $\mathcal{F}^{op}:\mathbb{R}^{n}\times 
\mathbb{R}^{m}\rightrightarrows \mathbb{R}^{n},$ given by 
\begin{equation}
\mathcal{F}^{op}\left( c,b\right) :=\arg \min \left\{ c^{\prime }x\mid
Ax\leq b\right\} .  \label{eq_Fop}
\end{equation}%
The parameterized linear system associated with (\ref{eq_LPproblem}) gives
rise to the \emph{feasible set mapping }$\mathcal{F}:\mathbb{R}%
^{m}\rightrightarrows \mathbb{R}^{n},$ which is defined as 
\begin{equation}
\mathcal{F}\left( b\right) :=\left\{ x\in \mathbb{R}^{n}\mid Ax\leq
b\right\} ,\text{ }b\in \mathbb{R}^{m}.  \label{eq_F(b)}
\end{equation}

At this moment let us introduce some notation. Given any metric space $%
\left( Z,d\right) ,$ the point-to-set distance from any element $x\in Z$ to
any subset $\Omega \subset Z$ is defined, as usual, by $d\left( x,\Omega
\right) :=\inf \left\{ d\left( x,z\right) :z\in Z\right\} ,$ with the
convention $\inf \emptyset :=+\infty ,$ so that $d\left( x,\emptyset \right)
=+\infty .$ Let $\mathcal{S}:Y\rightrightarrows X$ be a \emph{set-valued
mapping} (also called \emph{multifunction}) between the metric spaces $Y$
and $X$, with distance functions denoted by $d_{Y}$ and $d_{X},$ respectively%
$.\,$The \emph{domain} and the \emph{graph} of $\mathcal{S}$ are,
respectively, $\limfunc{dom}\mathcal{S}=\left\{ y\in Y\mid \mathcal{S}\left(
y\right) \neq \emptyset \right\} $ and $\limfunc{gph}\mathcal{S}=\left\{
\left( y,x\right) \in Y\times X\mid x\in \mathcal{S}\left( y\right) \right\}
;$ the \emph{inverse multifunction} $\mathcal{S}^{-1}:X\rightrightarrows Y$
is given by $y\in \mathcal{S}^{-1}\left( x\right) \Leftrightarrow x\in 
\mathcal{S}\left( y\right) $. We say that $\mathcal{S}$ is \emph{Lipschitz
continuous} \emph{on its domain }if\emph{\ }there exists a constant $\kappa
\geq 0$ such that%
\begin{equation}
d_{X}\left( x_{1},\mathcal{S}\left( y_{2}\right) \right) \leq \kappa
d_{Y}\left( y_{1},y_{2}\right) \text{ for all }y_{1},y_{2}\in \mathrm{dom}%
\mathcal{S}\text{ and all }x_{1}\in \mathcal{S}\left( y_{1}\right) .
\label{eq_Hoffman_constant}
\end{equation}%
Along this work, any constant $\kappa \geq 0$ verifying (\ref%
{eq_Hoffman_constant}) is referred to as a \emph{Lipschitz constant }of $%
\mathcal{S}.$ The infimum of all Lipschitz constants of $\mathcal{S}$ is the
so-called \emph{Hoffman constant} of $\mathcal{S},$ denoted by $\Hof\mathcal{%
S}.$ It is clear that $\Hof\mathcal{S=+\infty }$ when $\mathcal{S}$ is not 
\emph{Lipschitz continuous} \emph{on its domain}. From the definition, one
can easily check that 
\begin{equation}
\Hof\mathcal{S}=\sup\limits_{\substack{ \left( y_{1},x_{1}\right) \in 
\limfunc{gph}\mathcal{S},  \\ y_{2}\in \mathrm{dom}\mathcal{S}}}~~\dfrac{%
d_{X}\left( x_{1},\mathcal{S}\left( y_{2}\right) \right) }{d_{Y}\left(
y_{1},y_{2}\right) }=\sup\limits_{\substack{ x\in X,\text{ }  \\ y\dom%
\mathcal{S}}}~~\dfrac{d_{X}(x,\mathcal{S}(y))}{d_{Y}\left( y,\mathcal{S}%
^{-1}(x)\right) },  \label{eq_Hof}
\end{equation}%
under the convention $\frac{0}{0}:=0.\,$

The terminology `\emph{Hoffman constant}' is named after the celebrated
Hoffman lemma, established in \cite{Hof52} (see also G\"{u}ler, Hoffman, and
Rothblum \cite[Theorem 1.1]{GHR95}), which entails the existence of a
Lipschitz constant for the feasible set mapping introduced in (\ref{eq_F(b)}%
) and, consequently, the finiteness of $\Hof\mathcal{F}.$ Furthermore, we
emphasize that $\Hof\mathcal{F}$ is closely related to condition measures
and we will explore this relationship more in detail in Subsection 6.1 by
referring to Ho and Tun\c{c}el \cite{HoTu02} as well as to Ekbatani, Natura,
and V\'{e}gh \cite{ENV22}.

Since the pioneer work \cite{Hof52}, the interest on this topic is evidenced
by the remarkable research on it; in fact, exact formulae for $\Hof\mathcal{F%
}$ (under different choices of distances in $\mathbb{R}^{m}$ and $\mathbb{R}%
^{n}$) can be traced out from Burke and Tseng \cite[Theorem 8]{BurkeTseng96}%
, Klatte and Thiere \cite[Theorem 2.7]{KlTh95}, Li \cite[Theorems 2.4 and 3.4%
]{Li93} and Pe\~{n}a, Vera, and Zuluaga \cite[Formula (3)]{PVZ20}, among
others; see \cite{CCP21} for the extension to the semi-infinite setting
(with infinitely many inequality constraints; see also the monograph \cite%
{libro}). Let us also cite other relevant contributions by Az\'{e} and
Corvellec \cite{AzCo02}, Bergthaller and Singer \cite{BeSi92}, Cook,
Gerards, Schrijver, and Tardos \cite{CoGeScrhi86}, Li \cite{Li94},
Mangasarian and Shiau \cite{ManShian87}, and Z\u{a}linescu \cite{Za03}.

The present work is focused on the Hoffman constant for the argmin mapping
under different kinds of perturbations: the so-called \emph{canonical
perturbations }(with both $c$ and $b$ being considered as parameters that
can be perturbed), $c$ perturbations ($c$ is the parameter and $b$ is fixed)
and right-hand-side perturbation ($b$ is the parameter and $c$ is fixed). We
anticipate that the most interesting scenario is that of right-hand-side
(RHS, for short) perturbations, as it is shown in Section 3. Indeed,
Propositions \ref{Prop Hof canonical} and \ref{Prop_b_fixed} in Section 3
state that $\Hof\mathcal{F}^{op}$ and the Hoffman constant of 
\begin{equation}
\mathcal{F}_{\overline{b}}^{op}:=\mathcal{F}^{op}\left( \cdot ,\overline{b}%
\right) ,  \label{eq_F_b}
\end{equation}%
for a fixed $\overline{b}\in \mathbb{R}^{m},$ are either $0$ (in very
particular cases) or $+\infty $ . Hence, the main part of this paper is
concerned with the optimal set mapping under RHS perturbations, $\mathcal{F}%
_{\overline{c}}^{op}:\mathbb{R}^{m}\rightrightarrows \mathbb{R}^{n}$ with $%
\overline{c}\in \mathbb{R}^{n}$ fixed, defined as 
\begin{equation}
\mathcal{F}_{\overline{c}}^{op}\left( b\right) :=\mathcal{F}^{op}\left( 
\overline{c},b\right) \text{, for any }b\in \mathbb{R}^{m}.
\label{eq_F_op_c_fixed}
\end{equation}%
To be precise, most of the work is concerned with the computation of $\Hof%
\mathcal{F}_{\overline{c}}^{op}$ through a \emph{point-based} formula; i.e.,
through some expression involving only the fixed (also called \emph{nominal}%
) data $\overline{c}$ and $A$ (recall that $b$ is the only parameter in this
context) and the norms being considered in $\mathbb{R}^{n}$ and $\mathbb{R}%
^{m}.\,$Observe that $\Hof\mathcal{F}_{\overline{c}}^{op}$ provides a \emph{%
global} measure of the stability of optimal solutions with respect to
parameter perturbations. In fact, roughly speaking, this Hoffman constant is
the `tightest global rate of variation' of the set of optimal solutions with
respect to perturbations of parameter $b$.

We can find in the literature upper estimates of this constant as the one
provided in \cite[Section 6]{HoTu02} and, to the authors' knowledge, the
current paper establishes in Theorem \ref{Th upper bound Hof} (Subsection
4.3) the first exact formula for the aimed $\Hof\mathcal{F}_{\overline{c}%
}^{op}$. Looking at this formula, one easily sees that $\Hof\mathcal{F}_{%
\overline{c}}^{op}<+\infty ,$ yielding the Lipschitz continuity of $\mathcal{%
F}_{\overline{c}}^{op}$ restricted to its domain. In fact, the fulfilment of
such a property is known since the 1980s; see e.g. \cite{CoGeScrhi86,
DoZo93, Kla87, Li93, Li94, Man81, ManShian87} (see also \cite[Sec. 4]{KlTh95}
for an extension to quadratic problems), some of them even providing
specific expressions for Lipschitz constants (upper estimates for $\Hof%
\mathcal{F}_{\overline{c}}^{op}$).

Hoffman constants of the feasible and the optimal set mapping play a
fundamental role in mathematical programming, in particular, regarding the
stability analysis and convergence properties of a variety of optimization
algorithms; regarding the feasible set, see, e.g., \cite[Sections 7 and 8]%
{ENV22}, as well as \cite{PVZ20} and references therein (e.g., \cite%
{LuoTseng, NeNesGli, PeRo19, ZHOUSO17}) and with respect to the optimal set,
observe that it is closely related to feasibility of linear complementarity
constraints via the equivalence between optimal solutions and the well-known
Karush-Kuhn-Tucker points for linear programs. Going further, Lipschitz type
properties are in the core of variational analysis and the reader is
addressed, among others, to the monographs of Dontchev and Rockafellar \cite%
{DoRo}, Ioffe \cite{Io17}, Klatte and Kummer \cite{KlKu02}, Mordukhovich 
\cite{mor06a} and Rockafellar and Wets \cite{rw} for comprehensive
developments in this field.

Now we discuss the structure of the paper and some notable aspects of its
methodology. Section \ref{sec: pre} provides some notation and preliminary
results used throughout the work, and comments more in detail on condition
measures, as pointed out above. We consider some local, semilocal and global
variational rates (the terminology is explained in Section 2). Section 3 is
devoted to compute the Hoffman constants of $\mathcal{F}^{op}$ and $\mathcal{%
F}_{\overline{b}}^{op}$ which are $+\infty $ except in very particular cases
where they are equal to zero$.$ Inspired by what happens to the feasible
set, Section 4 tackles the computation of the Hoffman constant of $\mathcal{F%
}_{\overline{c}}^{op}$ by answering the following question: Is it possible
to relate the global stability measure $\Hof\mathcal{F}_{\overline{c}}^{op}$
with previously known \emph{local }ones? Paper \cite{CCP21} gives a positive
answer for the feasible set mapping as $\Hof\mathcal{F}$ is shown to
coincide with the supremum of the \emph{calmness moduli} of $\mathcal{F}$ at
all $b\in \dom\mathcal{F}$ and all $x\in \mathcal{F}\left( b\right) $ (see
Section 2 for the formal definition).\thinspace\ In that work, the convexity
of the graph of $\mathcal{F}$ is a key fact. Regarding our optimal set
mapping $\mathcal{F}_{\overline{c}}^{op},$ its graph is not convex in
general and this paper develops \emph{ad hoc }techniques to answer the
previous question. These techniques are developed in three subsections:
Subsection 4.1 analyzes the particular structure of the set-valued mapping $%
\mathcal{F}_{\overline{c}}^{op}$ as union of mappings with convex graphs
which are connected in a special manner. Subsection 4.2 formalizes this
particular structure giving rise to the concept of \emph{well-connected
piecewise convex mapping. }This concept is intended to isolate the main
features behind the possibility of relating the (global) Hoffman constant
with the (local) calmness moduli. The results of Subsection 4.2 apply in the
following subsection to compute the aimed $\Hof\mathcal{F}_{\overline{c}%
}^{op}$ through the announced calmness moduli; a remarkable fact is that
formulas for these calmness moduli already exist (see Theorem \ref{theo: clm
feasible}). Section 5 introduces the concept of break steps and provides a
constructive procedure to connect any two RHS parameters $b,$ $\overline{b}%
\in \limfunc{dom}\mathcal{F}_{\overline{c}}^{op}$, yielding, in Proposition %
\ref{Theorem_sect4}, an upper estimation of the variation of the optimal set
between those two parameters. We finish the paper with a section of
conclusions, perspectives, and examples. This last section includes examples
and arguments showing that previous constant introduced \cite[Section 6]%
{HoTu02} and \cite{Li93} can be strictly greater than $\Hof\mathcal{F}_{%
\overline{c}}^{op}.$

\section{Notation, preliminaries and first results}

\label{sec: pre}

We introduce some more notation used hereafter. Given $S\subset \mathbb{R}%
^{p}$, $p\in \mathbb{N}$, we denote by $\mathrm{int}S,$ $\mathrm{cl}S,$ $%
\mathrm{conv}S,$ $\mathrm{cone}S,$ and $\mathrm{span}S$ the interior, the
closure, the convex hull, the conical convex hull, and the linear subspace
of $\mathbb{R}^{n}$ spanned by $S,$ respectively, with the convention $%
\mathrm{conv}\emptyset =\emptyset $ and $\mathrm{cone}\emptyset =\spann%
\emptyset =\left\{ 0_{p}\right\} $ (the zero vector of $\mathbb{R}^{p}$).
Provided that $S$ is convex, $\mathrm{extr}S$ stands for the set of extreme
points of $S$.

Now we recall the \emph{semilocal }and \emph{local }Lipschitz type
properties referenced in the paper for a multifunction $\mathcal{S}%
:Y\rightrightarrows X$ between metric spaces, $Y$ and $X,$ with distance
functions being denoted by $d_{Y}$ and $d_{X},$ respectively$.$ With respect
to semilocal properties, this paper focuses on the \emph{Lipschitz upper
semicontinuity }of $\mathcal{S}$ at $\overline{y}\in \dom\mathcal{S},$ which
is defined as the existence of a neighborhood $V$ of $\overline{y}$ along
with a constant $\kappa \geq 0$ such that%
\begin{equation}
d_{X}(x,\mathcal{S}(\overline{y}))\leq \kappa d_{Y}\left( y,\overline{y}%
\right) \text{ for all }y\in V\text{ and all }x\in \mathcal{S}\left(
y\right) .  \label{eq_008}
\end{equation}%
Here the terminology `semilocal' means that $y$ varies around a fixed
element, $\overline{y},$ while $x$ runs over the whole image set $\mathcal{S}%
\left( y\right) $. The infimum of constants $\kappa ,$ for some associated
neighborhoods, appearing in (\ref{eq_008}) is the \emph{Lipschitz upper
semicontinuity modulus }of\emph{\ }$\mathcal{S}$ at $\overline{y}.$ It is
denoted by $\Lipusc\mathcal{S}(\overline{y})$ and \cite[Proposition 2]{CCP21}
establishes the following expression for it%
\begin{equation}
\Lipusc\mathcal{S}(\overline{y})=\limsup\limits_{y\rightarrow \overline{y}%
}\left( \sup\limits_{x\in \mathcal{S}(y)}\dfrac{d_{X}(x,\mathcal{S}(%
\overline{y}))}{d_{Y}\left( y,\overline{y}\right) }\right) ,\overline{y}\in %
\dom\mathcal{S}.  \label{eq_Lipusc}
\end{equation}%
The concept of $\limsup_{z\rightarrow \overline{z}}$ (with $z$ and $%
\overline{z}$ in a metric space) is standard, and in some proofs it is
useful to write it as the supremum (maximum, in fact) of all possible
sequential upper limits for all possible sequences $\left\{ z_{r}\right\}
_{r\in \mathbb{N}}$ converging to $\overline{z}$ as $r\rightarrow \infty $.
Recall the convention $\frac{0}{0}:=0$.

We also deal with the \emph{calmness} property, which is a local measure as
it considers solutions near a given solution $\overline{x}$ and parameters
in a neighborhood of the nominal one $\overline{y}$. Specifically, $\mathcal{%
S}$ is said to be calm at $\left( \overline{y},\overline{x}\right) \in 
\mathrm{gph}\mathcal{S}$ (the graph of $\mathcal{S}$) if there exist a
constant $\kappa \geq 0$ and a neighborhood of $\left( \overline{y},%
\overline{x}\right) ,$ $V\times U,$ such that 
\begin{equation}
d_{X}(x,\mathcal{S}(\overline{y}))\leq \kappa d_{Y}\left( y,\overline{y}%
\right) \text{ for all }x\in \mathcal{S}(y)\cap U\text{ and all }y\in V,
\label{eq_calmM}
\end{equation}%
which is known to be equivalent to the \emph{metric subregularity }(cf. \cite%
[Theorem 3H.3 and Exercise 3H.4]{DoRo}) of $\mathcal{S}^{-1}$ at $\left( 
\overline{x},\overline{y}\right) ,$ which reads as the existence of $\kappa
\geq 0$ and a (possibly smaller) neighborhood $U$ of $\overline{x}$ such
that 
\begin{equation}
d_{X}(x,\mathcal{S}(\overline{y}))\leq \kappa d_{Y}\left( \overline{y},%
\mathcal{S}^{-1}(x)\right) \text{ for all }x\in U\text{ }.
\label{eq_met_subreg}
\end{equation}%
Provided that $\kappa >0,$ from \cite[Theorem 3H.3 and Exercise 3H.4]{DoRo}
we know that $\mathcal{S}$ is calm at $\left( \overline{y},\overline{x}%
\right) $ with associated calmness constant $\kappa $ if and only if $%
\mathcal{S}^{-1}$ is metrically subregular at $\left( \overline{x},\overline{%
y}\right) $ with constant $\kappa .$ Therefore, the infima of constants $%
\kappa ,$ for some associated neighborhoods, appearing in (\ref{eq_calmM})
and (\ref{eq_met_subreg}) coincide. This infimum is the so-called \emph{%
calmness modulus} of $\mathcal{S}$ at $\left( \overline{y},\overline{x}%
\right) ,$ denoted by $\clm\mathcal{S}\left( \overline{y},\overline{x}%
\right) .$ Directly from the definition, for any $\left( \overline{y},%
\overline{x}\right) \in \gph\mathcal{S},$ we have

\begin{equation*}
\clm\mathcal{S}\left( \overline{y},\overline{x}\right) =\limsup\limits 
_{\substack{ \left( x,y\right) \rightarrow \left( \overline{x},\overline{y}%
\right)  \\ \left( x,y\right) \in \func{gph}\mathcal{S}}}\dfrac{d_{X}(x,%
\mathcal{S}(\overline{y}))}{d_{Y}\left( y,\overline{y}\right) }%
=\limsup\limits_{x\rightarrow \overline{x}}\dfrac{d_{X}(x,\mathcal{S}(%
\overline{y}))}{d_{Y}\left( \overline{y},\mathcal{S}^{-1}(x)\right) }.
\end{equation*}

It is clear from the definitions that 
\begin{equation}
\Hof\mathcal{S\geq }\sup\limits_{y\in \dom\mathcal{S}}\Lipusc\mathcal{S}%
\left( y\right) \geq \sup_{\left( y,x\right) \in \mathrm{gph}\mathcal{S}}\clm%
\mathcal{S}\left( y,x\right) .  \label{eq_inequalities chain}
\end{equation}

The following theorem, which can be found in \cite[Theorem 4]{CCP21},
provides a sufficient condition to get equalities in (\ref{eq_inequalities
chain}); Corollary \ref{Cor_Hofmann_Lipusc} and Theorem \ref{Th upper bound
Hof} (in Section 4) constitute the version of this result adapted to
well-connected piecewise convex mappings in finite dimensions and for the
particular case of $\mathcal{F}_{\overline{c}}^{op},$ respectively.

\begin{theo}
\label{The_calmness_Hoffman_convex}Let $\mathcal{S}:Y\rightrightarrows X,$
with $Y$ being a normed space and $X$ being a reflexive Banach space. Assume
that $\gph\mathcal{S}$ is a nonempty convex set and that, for all $y\in \dom%
\mathcal{S},$ $\mathcal{S}\left( y\right) $ is a closed subset of $X$. Then 
\begin{equation}
\Hof\mathcal{S}=\sup_{y\in \dom\mathcal{S}}\Lipusc\mathcal{S}\left( y\right)
=\sup_{\left( y,x\right) \in \gph\mathcal{S}}\clm\mathcal{S}\left(
y,x\right) .  \label{eq_0001}
\end{equation}
\end{theo}

Clearly, $\gph\mathcal{F}$ is a convex set and $\mathcal{F}$ is
closed-valued (in fact, $\gph\mathcal{F}$ is a convex polyhedral set) and,
hence, the previous theorem applies for $\mathcal{S=F}.$ Regarding the
argmin mapping, $\gph\mathcal{F}^{op}$ is not convex and $\gph\mathcal{F}_{%
\overline{c}}^{op}$ is also nonconvex in general (see, e.g., \cite[Example
1.1]{CCP21Lip}). This underlies the fact that the analysis of $\Hof\mathcal{F%
}^{op}\mathcal{\ }$and $\Hof\mathcal{F}_{\overline{c}}^{op}$ does not rely
on Theorem \ref{The_calmness_Hoffman_convex}. As announced in Section 1,
this analysis constitutes the initial goal of the current paper and it is
developed in Section 4 applying the results of Section 3 about
well-connected piecewise convex mappings.\emph{\ }

For completeness, the rest of this section is devoted to provide some
background on the calmness and Lipschitz upper semicontinuity moduli for
multifunction $\mathcal{F}^{op}.$ Specifically, Theorems \ref{theo: clm
feasible} and \ref{theo: Lip extreme points} provide point-based formulae
(only depending on the nominal parameter and point, and not elements in a
neighborhood) for these two constants. Although $\gph\mathcal{F}^{op}$ is
not convex, $\mathcal{F}^{op}$ still satisfies a certain local directional
convexity property, which turns out to be crucial for obtaining Theorem \ref%
{theo: Lip extreme points} (see \cite[Theorem 5]{CCP21Lip}).

Let us introduce some more notation and fix the topology of the involved
spaces. The space of variables, $\mathbb{R}^{n},$ is endowed with an
arbitrary norm $\Vert \cdot \Vert $, with dual norm $\Vert \cdot \Vert
_{\ast }$, whereas the RHS parameter space, $\mathbb{R}^{m}$, is endowed
with the maximum (Chebyshev) norm $\Vert \cdot \Vert _{\infty }$. The
distances associated with $\Vert \cdot \Vert $ and $\Vert \cdot \Vert _{\ast
}$ in $\mathbb{R}^{n}$ are denoted by $d$ and $d_{\ast },$ respectively;
whereas the distance in $\mathbb{R}^{m}$ associated with the norm $\Vert
\cdot \Vert _{\infty }$ is denoted as $d_{\infty }.$ The full parameter
space, $\mathbb{R}^{n}\times \mathbb{R}^{m}$, is endowed with the norm 
\begin{equation*}
\Vert \left( c,b\right) \Vert :=\max \left\{ \Vert c\Vert _{\ast },\Vert
b\Vert _{\infty }\right\} ,
\end{equation*}
since $c$ is regarded as the linear functional $x\mapsto c^{\prime }x.$ In
our analysis it will be convenient to name explicitly the rows of $A;$ so
that $a_{t}^{\prime }$ (transpose of $a_{t}\in \mathbb{R}^{n}$) denotes the $%
t$-th row of $A$. Hence, the constraint system of problem (\ref{eq_LPproblem}%
) reads as 
\begin{equation}
a_{t}^{\prime }x\leq b_{t},~~t\in T:=\left\{ 1,...,m\right\} .
\label{system ab}
\end{equation}%
From now on we refer to the set of active indices at $x\in \mathcal{F}\left(
b\right) ,$ defined as 
\begin{equation*}
T_{b}\left( x\right) :=\left\{ t\in T\mid a_{t}^{\prime }x=b_{t}\right\} .
\end{equation*}

\begin{defn}
\label{def_KKT_subsets}$\left( i\right) $ Given $\left( c,b\right) \in \dom%
\mathcal{F}^{op}$ and $x\in \mathcal{F}^{op}\left( c,b\right) ,$ the family
of \emph{minimal KKT index subsets} at $\left( \left( c,b\right) ,x\right) $
denoted by $\mathcal{M}_{c,b}\left( x\right) $ --and introduced in \cite%
{CHLP16}--, is defined as the collection of all $D\subset T_{b}\left(
x\right) $ such that $D$ is minimal (with respect to the inclusion order)
among those satisfying $-c\in \cone\left\{ a_{t},\,t\in D\right\} .$

$\left( ii\right) $ It has been proved that $\mathcal{M}_{c,b}\left(
x\right) $ does not depend on $x\in \mathcal{F}^{op}\left( c,b\right) $ (cf. 
\cite[Remark 2]{GCPT18}), so that it is just referred to as the family of 
\emph{minimal KKT index subsets} at $\left( c,b\right) $, denoted as $%
\mathcal{M}_{c,b}.$
\end{defn}

\begin{lem}
\emph{\cite[Lemma 3.2]{CCP21Lip}}\label{lem: tech} Let $(\overline{c},%
\overline{b})\in \dom\mathcal{F}^{op}$. Then there exists $\varepsilon >0$
such that for every $b\in \dom\mathcal{F}$ with $\Vert b-\overline{b}\Vert
_{\infty }\leq \varepsilon $ we have $\mathcal{M}_{\overline{c},b}\subset 
\mathcal{M}_{\overline{c},\overline{b}}$.
\end{lem}

For any $D\subset T$ we consider the mapping $\mathcal{L}_{D}:\mathbb{R}%
^{m}\times \mathbb{R}^{D}\rightrightarrows \mathbb{R}^{n}$ given by 
\begin{equation}
\mathcal{L}_{D}\left( b,d\right) :=\left\{ x\in \mathbb{R}^{n}\mid
a_{t}^{\prime }x\leq b_{t},\,t\in T;\;-a_{t}^{\prime }x\leq d_{t},\,t\in
D\right\} .  \label{eq_003}
\end{equation}%
Observe that $\mathcal{L}_{D}$ is nothing else but the feasible set mapping
associated with an extension of the constraint system of (\ref{eq_LPproblem}%
) in order to force that inequalities indexed by $D$ are held as equalities
at the nominal parameter. In what follows $\overline{b}_{D}$ means $\left( 
\overline{b}_{t}\right) _{t\in D}.$

\begin{prop}
\emph{\cite[Proposition 4.1]{CHLP16}}\label{Prop LD=Fop} Let $(\overline{c},%
\overline{b})\in \dom\mathcal{F}^{op}.$ Then 
\begin{equation*}
\mathcal{L}_{D}\left( \overline{b},-\overline{b}_{D}\right) =\mathcal{F}%
^{op}(\overline{c},\overline{b})\text{ for all }D\in \mathcal{M}_{\overline{c%
},\overline{b}}.
\end{equation*}
\end{prop}

The next result provides three different expressions for the calmness
modulus of the optimal set mapping $\mathcal{F}^{op}$ at $\left( \left( 
\overline{c},\overline{b}\right) ,\overline{x}\right) \in \gph\mathcal{F}%
^{op}$. The first two ones come directly from \cite[Corollary 4.1]{CHLP16};
in the second expression each $\clm\mathcal{L}_{D}\left( \left( \overline{b}%
,-\overline{b}_{D}\right) ,\overline{x}\right) $ can be computed through the
concept of \emph{end set} of a convex set $C\subset \mathbb{R}^{n}$,
introduced in \cite{Hu05}\emph{\ }(see also \cite{ZhNg04}){\Huge \ }and
defined as 
\begin{equation*}
\mathrm{end\,}C:=\left\{ u\in \mathrm{cl\,}C\mid \nexists \mu >1\text{ such
that }\mu u\in \mathrm{cl\,}C\right\} .
\end{equation*}%
The third expression for $\clm\mathcal{F}^{op}\left( \left( \overline{c},%
\overline{b}\right) ,\overline{x}\right) $ below can be seen as a
geometrical interpretation of the formula given in \cite[Theorem 4]{CLPT14}
for the calmness modulus of a feasible set mapping.

\begin{theo}
\emph{\cite[Corollary 4.1]{CHLP16}, \cite[Theorem 4]{CLPT14}} \label{theo:
clm feasible} Let $\left( \left( \overline{c},\overline{b}\right) ,\overline{%
x}\right) \in \gph\mathcal{F}^{op}$. Then 
\begin{multline*}
\clm\mathcal{F}^{op}\left( \left( \overline{c},\overline{b}\right) ,%
\overline{x}\right) =\clm\mathcal{F}_{\overline{c}}^{op}\left( \overline{b},%
\overline{x}\right) =\max_{D\in \mathcal{M}_{\overline{c},\overline{b}}}\clm%
\mathcal{L}_{D}\left( \left( \overline{b},-\overline{b}_{D}\right) ,%
\overline{x}\right) \\
=\max_{D\in \mathcal{M}_{\overline{c},\overline{b}}}\left[ d_{\ast }\left(
0_{n},\mathrm{end\,conv}\left\{ a_{t},~t\in T_{\overline{b}}\left( \overline{%
x}\right) ;-a_{t},~t\in D\right\} \right) \right] ^{-1}.
\end{multline*}
\end{theo}

Let us introduce a natural extension for the set of extreme points: 
\begin{align*}
\mathcal{E}\left( b\right) & :=\extr\left( \mathcal{F}\left( b\right) \cap %
\spann\left\{ a_{t},\,t\in T\right\} \right) ,\;b\in \dom\mathcal{F}, \\
\mathcal{E}^{op}\left( c,b\right) & :=\extr\left( \mathcal{F}^{op}\left(
c,b\right) \cap \spann\left\{ a_{t},\,t\in T\right\} \right) ,\;\left(
c,b\right) \in \dom\mathcal{F}^{op}.
\end{align*}%
The reader is addressed to \cite[p. 142]{Li94} and \cite[Section 2.2]{GCPT18}
for details about these constructions. Note that $\mathcal{E}^{op}\left(
c,b\right) =\mathcal{F}^{op}\left( c,b\right) \cap \mathcal{E}\left(
b\right) $ for $\left( c,b\right) \in \dom\mathcal{F}^{op}$ and it is a
nonempty finite set.

\begin{theo}
\emph{\cite[Corollary 4.1, Proposition 4.2, and Theorem 4.2]{CCP21Lip}} %
\label{theo: Lip extreme points} Let $\left( \overline{c},\overline{b}%
\right) \in \dom\mathcal{F}^{op}$, then 
\begin{eqnarray*}
\Lipusc\mathcal{F}^{op}\left( \overline{c},\overline{b}\right) &=&\Lipusc%
\mathcal{F}_{\overline{c}}^{op}\left( \overline{b}\right) \\
&=&\sup_{x\in \mathcal{F}^{op}\left( \overline{c},\overline{b}\right) }\clm%
\mathcal{F}^{op}\left( \left( \overline{c},\overline{b}\right) ,x\right)
=\max_{x\in \mathcal{E}^{op}\left( \overline{c},\overline{b}\right) }\clm%
\mathcal{F}^{op}\left( \left( \overline{c},\overline{b}\right) ,x\right) .
\end{eqnarray*}
\end{theo}

\section{Hoffman constant under canonical and c perturbations}

This section deals with the optimal set mappings under canonical and $c$
perturbations, $\mathcal{F}^{op}$ and $\mathcal{F}_{\overline{b}}^{op},$
introduced in (\ref{eq_Fop}) and (\ref{eq_F_b}), respectively. Indeed, it is
oriented to show that $\Hof\mathcal{F}^{op}$ and $\Hof\mathcal{F}_{\overline{%
b}}^{op}$ are infinite unless we are placed very particular cases, as
formalized in Propositions \ref{Prop Hof canonical} and \ref{Prop_b_fixed}.

Recall that, in ordinary (finite) linear programming, optimality is
equivalent to primal-dual consistency. In other words, 
\begin{equation}
\dom\mathcal{F}^{op}=\left( -\cone\left\{ a_{t},\,t\in T\right\} \right)
\times \dom\mathcal{F}.  \label{eq_domFop}
\end{equation}%
Clearly, $\gph\mathcal{F},$ $\dom\mathcal{F},$ and $\dom\mathcal{F}^{op}$
are convex, while $\gph\mathcal{F}^{op}$ is not convex. In the next
propositions $`0_{m\times n}$' denotes the matrix in $\mathbb{R}^{m\times n}$
with all its coordinates equal to zero.

\begin{prop}
\label{Prop Hof canonical}We have 
\begin{equation*}
\Hof\mathcal{F}^{op}=\left\{ 
\begin{array}{l}
0\text{ if }A=0_{m\times n}, \\ 
+\infty \text{ otherwise.}%
\end{array}%
\right.
\end{equation*}
\end{prop}

\begin{dem}
In the case when $A=0_{m\times n},$ one trivially has 
\begin{equation*}
\mathcal{F}^{op}\left( c,b\right) =\left\{ 
\begin{array}{l}
\mathbb{R}^{n}\text{ if }\left( c,b\right) \in \{0_{n}\}\times \mathbb{R}%
_{+}^{m}, \\ 
\emptyset \text{ otherwise.}%
\end{array}%
\right.
\end{equation*}%
So, it is clear that $\Hof\mathcal{F}^{op}=0$.

Assume now that $\left\{ a_{t},\text{ }t\in T\right\} \neq \{0_{n}\}.$ Take $%
\overline{x}\in \mathbb{R}^{n}$ and define $\overline{b}_{t}=a_{t}^{\prime }%
\overline{x}+1$ for all $t\in T.$ Hence $\overline{x}\in \mathrm{int}%
\mathcal{F}\left( \overline{b}\right) $ and it cannot be an optimal solution
for any $\left( c,\overline{b}\right) \in \dom\mathcal{F}^{op}$ with $c\neq
0_{n}.$ Fix any $\overline{c}\in -\cone\left\{ a_{t},\,t\in T\right\}
\setminus \{0_{n}\}$. We have $d\left( \overline{x},\mathcal{F}^{op}\left( 
\overline{c},\overline{b}\right) \right) =d\left( \overline{x},\mathcal{F}%
^{op}\left( \frac{1}{r}\overline{c},\overline{b}\right) \right) >0$ for all $%
r\in \mathbb{N}.$

On the other hand, it is clear that $\overline{x}\in \mathcal{F}^{op}\left(
0_{n},\overline{b}\right) .$ Hence%
\begin{equation*}
\Hof\mathcal{F}^{op}\geq \lim_{r\rightarrow +\infty }\frac{d\left( \overline{%
x},\mathcal{F}^{op}\left( \frac{1}{r}\overline{c},\overline{b}\right)
\right) }{\left\Vert \left( 0_{n},\overline{b}\right) -\left( \frac{1}{r}%
\overline{c},\overline{b}\right) \right\Vert }=+\infty .
\end{equation*}
\end{dem}

\begin{prop}
\label{Prop_b_fixed}Let $\overline{b}\in \dom\mathcal{F}.$ We have 
\begin{equation*}
\Hof\mathcal{F}_{\overline{b}}^{op}=\left\{ 
\begin{array}{l}
0\text{ if }\mathcal{F}\left( \overline{b}\right) \text{ is an affine
variety of }\mathbb{R}^{n}, \\ 
+\infty \text{ otherwise.}%
\end{array}%
\right.
\end{equation*}
\end{prop}

\begin{dem}
Let $L$ be the lineality space of $\mathcal{F}\left( \overline{b}\right) $;
i.e., 
\begin{equation*}
L=\left( -0^{+}\mathcal{F}\left( \overline{b}\right) \right) \cap 0^{+}%
\mathcal{F}\left( \overline{b}\right) ,
\end{equation*}%
where $0^{+}\mathcal{F}\left( \overline{b}\right) $ is the recession cone of 
$\mathcal{F}\left( \overline{b}\right) $; see \cite[p. 65]{Rock70}, which
also shows that $\mathcal{F}\left( \overline{b}\right) $ can be expressed as
the direct sum 
\begin{equation*}
\mathcal{F}\left( \overline{b}\right) =L+\left( \mathcal{F}\left( \overline{b%
}\right) \cap L^{\bot }\right) ,
\end{equation*}%
with $L^{\bot }$ denoting the orthogonal complement of $L,$ which in our
case writes as $L^{\bot }=\limfunc{span}\left\{ a_{t},~t\in T\right\} .$ The
set of extreme points of $\mathcal{F}\left( \overline{b}\right) \cap L^{\bot
}$ is $\mathcal{E}\left( \overline{b}\right) ,$ introduced at the end of
Section 2, which is a finite nonempty set. Take any $\overline{x}\in 
\mathcal{E}\left( \overline{b}\right) .$ Then, $\overline{x}$ is an exposed
point of $\mathcal{F}\left( \overline{b}\right) \cap L^{\bot }$ according to 
\cite[Theorem 18.6]{Rock70}, which translates into the existence of $u\in
L^{\bot }\backslash \{0_{n}\}$ such that $u^{\prime }x\geq u^{\prime }%
\overline{x}$ for all $x\in \mathcal{F}\left( \overline{b}\right) \cap
L^{\bot }$ and 
\begin{equation*}
\mathcal{F}\left( \overline{b}\right) \cap L^{\bot }\cap \left\{ x\in 
\mathbb{R}^{n}:u^{\prime }x=u^{\prime }\overline{x}\right\} =\left\{ 
\overline{x}\right\} .
\end{equation*}%
Now we distinguish two cases:

\emph{Case 1}: $\mathcal{F}\left( \overline{b}\right) $ is an affine variety
of $\mathbb{R}^{n}.$ This is equivalent to $\mathcal{F}\left( \overline{b}%
\right) =\overline{x}+L;$ in other words, $\mathcal{F}\left( \overline{b}%
\right) \cap L^{\bot }=\left\{ \overline{x}\right\} .$ Then $c\in \limfunc{%
dom}\mathcal{F}_{\overline{b}}^{op}\Leftrightarrow c\in L^{\bot
}\Leftrightarrow \mathcal{F}_{\overline{b}}^{op}\left( c\right) =\mathcal{F}%
\left( \overline{b}\right) .$ Hence (\ref{eq_Hof}) clearly entails $\Hof%
\mathcal{F}_{\overline{b}}^{op}=0.$

\emph{Case 2}: $\mathcal{F}\left( \overline{b}\right) $ is not an affine
variety of $\mathbb{R}^{n}.$ In such a case take any $\widetilde{x}\in
\left( \mathcal{F}\left( \overline{b}\right) \cap L^{\bot }\right)
\backslash \{\overline{x}\}.$ Clearly $\widetilde{x}\in \mathcal{F}_{%
\overline{b}}^{op}\left( 0_{n}\right) =\mathcal{F}\left( \overline{b}\right)
.$ Then 
\begin{equation*}
\Hof\mathcal{F}_{\overline{b}}^{op}\geq \lim_{r\rightarrow \infty }\frac{%
d\left( \widetilde{x},\overline{x}+L\right) }{\left\Vert 0_{n}-\frac{1}{r}%
u\right\Vert _{\ast }}=\lim_{r\rightarrow \infty }r\frac{d\left( \widetilde{x%
},\overline{x}+L\right) }{\left\Vert u\right\Vert _{\ast }}=+\infty .
\end{equation*}
\end{dem}

\section{Hoffman constant under RHS perturbations}

This section is devoted to the computation of the Hoffman constant of $%
\mathcal{F}_{\overline{c}}^{op}$ through an implementable formula only
involving the fixed data $\overline{c}$ and matrix $A.$ The section is
divided into three subsections. The first one is intended to point out the
particular structure of $\mathcal{F}_{\overline{c}}^{op}$ which is behind
the fact that $\Hof\mathcal{F}_{\overline{c}}^{op}$ can be expressed as the
supremum of calmness moduli. This particular structure motivates the concept
of \emph{well-connected piecewise convex }mapping. The second subsection
analyzes the Hoffman constant for such a kind of set-valued mappings, which
allow us to deduce the aimed expression for $\Hof\mathcal{F}_{\overline{c}%
}^{op}$ from a more general perspective.

Along this section we assume that 
\begin{equation}
-\overline{c}\in \cone\left\{ a_{t},\,t\in T\right\} ,  \label{eq_222}
\end{equation}%
which entails, according to well-known arguments in ordinary (finite) linear
programming, that 
\begin{equation*}
\dom\mathcal{F}_{\overline{c}}^{op}=\dom\mathcal{F},
\end{equation*}%
hence $\dom\mathcal{F}_{\overline{c}}^{op}$ is a closed and convex set in $%
\mathbb{R}^{m}.$ From (\ref{eq_222}), the following family of subsets of
indices is nonempty:%
\begin{equation*}
\mathcal{K}_{\overline{c}}:=\left\{ D\subset T\left\vert -\overline{c}\in %
\cone\left\{ a_{t},\,t\in D\right\} \right. \right\} .
\end{equation*}%
From now on, given any subset $J\subset T,$ $A_{J}$ is the submatrix of $%
A\in \mathbb{R}^{m\times n}$ formed by the rows of $A$ indexed by $J$ and,
for any $b\in \mathbb{R}^{m},$ $b_{J}$ $\in \mathbb{R}^{J}$ is formed by the
coordinates of $b$ indexed by $J$.

The well-known Karush-Kuhn-Tucker (KKT in brief) optimality conditions
characterize those elements $\left( \overline{b},\overline{x}\right) \in \gph%
\mathcal{F}_{\overline{c}}^{op}$ by the existence of some $D\in \mathcal{K}_{%
\overline{c}}$ such that $A_{D}~\overline{x}=\overline{b}_{D},$ $%
A_{T\backslash D}~\overline{x}\leq \overline{b}_{T\backslash D};\;$in other
words, if, for each $D\in \mathcal{K}_{\overline{c}},$ we consider the
set-valued mapping $\mathcal{P}_{D}:\mathbb{R}^{m}\rightrightarrows \mathbb{R%
}^{n}$ assigning to each $b\in \mathbb{R}^{m}$ the convex polyhedral set 
\begin{equation}
\mathcal{P}_{D}\left( b\right) :=\left\{ x\in \mathbb{R}^{n}:A_{T\backslash
D}~x\leq b_{T\backslash D};\;A_{D}~x=b_{D}\right\} ,  \label{eq_004}
\end{equation}%
then, the KKT conditions read as 
\begin{equation}
\gph\mathcal{F}_{\overline{c}}^{op}=\tbigcup\limits_{D\in \mathcal{K}_{%
\overline{c}}}\gph\mathcal{P}_{D}.  \label{eq_KKT}
\end{equation}%
Observe that in general $\gph\mathcal{F}_{\overline{c}}^{op}$ is not convex,
but a finite union of convex polyhedral sets. It is clear that $\gph\mathcal{%
P}_{D}$ is a convex polyhedral (hence, closed) set and that $\dom\mathcal{P}%
_{D}$ is also a convex polyhedral set as the projection of $\gph\mathcal{P}%
_{D}\subset \mathbb{R}^{m}\times \mathbb{R}^{n}$ on $\mathbb{R}^{m}$.

\subsection{Piecewise structure of the argmin mapping under RHS perturbations%
}

This subsection provides a representation of $\mathcal{F}_{\overline{c}%
}^{op} $ as union of set-valued mappings coming from (\ref{eq_KKT}).

Given a family of set-valued mappings $\{\mathcal{S}_{i}:\mathbb{R}%
^{m}\rightrightarrows \mathbb{R}^{n},\ i\in I\},$ for an arbitrary set $I,$
the union mapping, $\tbigcup\nolimits_{i\in I}\mathcal{S}_{i}:\mathbb{R}%
^{m}\rightrightarrows \mathbb{R}^{n},$ is defined by 
\begin{equation*}
\left( \tbigcup\limits_{i\in I}\mathcal{S}_{i}\right) \left( b\right)
:=\tbigcup\limits_{i\in I}\mathcal{S}_{i}\left( b\right) ;
\end{equation*}%
it is clear that $\gph\left( \tbigcup\limits_{i\in I}\mathcal{S}_{i}\right)
=\tbigcup\limits_{i\in I}\gph\mathcal{S}_{i}$ and $\dom\left(
\tbigcup\limits_{i\in I}\mathcal{S}_{i}\right) =\tbigcup\limits_{i\in I}\dom%
\mathcal{S}_{i}.$

From (\ref{eq_KKT})$,$ one immediately has 
\begin{equation}
\mathcal{F}_{\overline{c}}^{op}=\tbigcup\limits_{D\in \mathcal{K}_{\overline{%
c}}}\mathcal{P}_{D}.  \label{eq_111}
\end{equation}

For any pair of subsets $D_{1},D_{2}\in \mathcal{K}_{\overline{c}},$ such
that $D_{1}\subset D_{2},$ it is clear that $\limfunc{dom}\mathcal{P}%
_{D_{2}}\subset \limfunc{dom}\mathcal{P}_{D_{1}},$ and $\mathcal{P}%
_{D_{2}}\left( b\right) \subset \mathcal{P}_{D_{1}}\left( b\right) $ for all 
$b\in \mathbb{R}^{m}.$ This motivates the fact of refining the previous
union in (\ref{eq_111}) by keeping those elements of $\mathcal{K}_{\overline{%
c}}$ which are minimal with respect to (w.r.t., in brief) the inclusion
order, as it is done in the next proposition. Formally, we consider the
following family of subsets of indices:%
\begin{equation}
\mathcal{M}_{\overline{c}}:=\left\{ D\in \mathcal{K}_{\overline{c}%
}\left\vert D\text{ is minimal w.r.t. the inclusion order}\right. \right\} .
\label{eq_002}
\end{equation}

\begin{prop}
\label{Prop_FOp_well_connected}We have:

$\left( i\right) \mathcal{F}_{\overline{c}}^{op}=\bigcup_{D\in \mathcal{M}_{%
\overline{c}}}\mathcal{P}_{D};$

$\left( ii\right) \mathcal{F}_{\overline{c}}^{op}\mid _{\dom\mathcal{P}_{D}}=%
\mathcal{P}_{D}$ for all $D\in \mathcal{M}_{\overline{c}}.$
\end{prop}

\begin{dem}
Condition $\left( i\right) $ comes directly from (\ref{eq_111}) and the
subsequent comments, while condition $\left( ii\right) $ can be derived from
Proposition \ref{Prop LD=Fop}. Specifically, if $D\in \mathcal{M}_{\overline{%
c}}$ and $b\in \dom\mathcal{P}_{D},$ then, it is clear that $D\in \mathcal{M}%
_{\overline{c},b}$ and, hence, Proposition \ref{Prop LD=Fop} yields 
\begin{equation*}
\mathcal{F}^{op}(\overline{c},b)=\mathcal{L}_{D}\left( b,-b_{D}\right) =%
\mathcal{P}_{D}\left( b\right) ,
\end{equation*}%
where the last equality comes directly from the definitions of the
set-valued mapping $\mathcal{P}_{D}$ and $\mathcal{L}_{D}$ (recall (\ref%
{eq_003})).
\end{dem}

\subsection{Hoffman constant for a class of piecewise convex mappings\protect%
\bigskip}

Motivated by the structure of $\mathcal{F}_{\overline{c}}^{op}$ stated in
Proposition \ref{Prop_FOp_well_connected}, we introduce the following
definition. This structure is behind the possibility of expressing the
Hoffman constant of a set-valued mapping by the supremum of the calmness
moduli at all points of its graph. Roughly speaking, a global variational
measure is given in terms of local ones.

\begin{defn}
\label{def_well_assembled}Let $I$ be a finite index set and, for each $i\in
I,$ consider a multifunction $\mathcal{S}_{i}:\mathbb{R}^{m}%
\rightrightarrows \mathbb{R}^{n}$ with $\gph\mathcal{S}_{i}$ and $\dom%
\mathcal{S}_{i}$ being closed convex sets in $\mathbb{R}^{m\times n}$ and $%
\mathbb{R}^{m},$ respectively. We say that $\mathcal{S}:=\tbigcup%
\nolimits_{i\in I}\mathcal{S}_{i}$ is a \emph{well-connected piecewise
convex }(\emph{wcpc}, for short)\emph{\ }mapping\emph{\ }if the following
properties are fulfilled:

$\left( i\right) $ $\dom\mathcal{S}$ $(=\tbigcup\nolimits_{i\in I}\dom%
\mathcal{S}_{i})$ is a convex set in $\mathbb{R}^{m};$

$\left( ii\right) $ $\mathcal{S}\mid _{\dom\mathcal{S}_{i}}=\mathcal{S}_{i},$
for all $i\in I$ (equivalently, $\mathcal{S}_{i}\left( b\right) =\mathcal{S}%
_{j}\left( b\right) $ whenever $b\in \dom\mathcal{S}_{i}\cap \dom\mathcal{S}%
_{j},$ $i,j\in I$).
\end{defn}

Throughout this subsection we consider $\mathbb{R}^{m}$ and $\mathbb{R}^{n}$
endowed with arbitrary norms, both denoted by $\left\Vert \cdot \right\Vert $
for simplicity and with $d$ standing for both associated distances. Given a
nonempty closed set $C\subset \mathbb{R}^{n}$, we denote by{}%
\begin{equation*}
P_{C}(x):=\mathrm{argmin}_{y}\{\left\Vert y-x\right\Vert \mid y\in C\}
\end{equation*}%
the set of best approximations (projections) of $x\in \mathbb{R}^{n}$ onto $%
C $. One easily checks (see \cite[Lemma 1]{CCP21}) that for any $\widetilde{x%
}\in $ $P_{C}(x)$ there holds%
\begin{equation}
\widetilde{x}\in P_{C}(\widetilde{x}+\mu \left( x-\widetilde{x}\right) ),%
\text{ for all }\mu \in [0,1].  \label{eq_00}
\end{equation}

The following lemma constitutes the version of \cite[Lemma 2]{CCP21} adapted
to our current context. In fact, this result holds for multifunctions with a
convex graph and a closed image set at the reference point (the closedness
of the whole graph is not required there). For completeness, we write a
sketch of the proof.

\begin{lem}
\label{Lem_ratio_dist}Let $\mathcal{S}=\tbigcup\nolimits_{i\in I}\mathcal{S}%
_{i}$ be a $wcpc$ mapping. Assume that $b,\widehat{b}\in \dom\mathcal{S}_{k}$
for some $k\in I.$ Take $x\in \mathcal{S}\left( b\right) $ and any $\widehat{%
x}\in P_{\mathcal{S}\left( \widehat{b}\right) }(x)$. Then%
\begin{equation*}
\frac{d\left( x,\mathcal{S}\left( \widehat{b}\right) \right) }{d\left( b,%
\widehat{b}\right) }\leq \clm\mathcal{S}\left( \widehat{b},\widehat{x}%
\right) .
\end{equation*}
\end{lem}

\begin{dem}
Since $\gph\mathcal{S}_{k}$ is a convex set, 
\begin{equation*}
\left( b_{\mu },x_{\mu }\right) :=(\widehat{b},\widehat{x})+\mu (\left(
b,x\right) -(\widehat{b},\widehat{x}))\in \gph\mathcal{S}_{k},\text{ for
each }\mu \in \left[ 0,1\right] .
\end{equation*}%
According to (\ref{eq_00}), $\widehat{x}\in P_{\mathcal{S}_{k}\left( 
\widehat{b}\right) }(x_{\mu })$, for each $\mu \in [ 0,1]$. Moreover, taking
into account that $\mathcal{S}_{k}\left( \widehat{b}\right) =\mathcal{S}%
\left( \widehat{b}\right) ,$ we have\textbf{\ } 
\begin{equation*}
\frac{d\left( x,\mathcal{S}\left( \widehat{b}\right) \right) }{d\left( b,%
\widehat{b}\right) }=\frac{\left\Vert x-\widehat{x}\right\Vert }{\left\Vert
b-\widehat{b}\right\Vert }=\frac{\left\Vert x_{\mu }-\widehat{x}\right\Vert 
}{\left\Vert b_{\mu }-\widehat{b}\right\Vert }=\frac{d\left( x_{\mu },%
\mathcal{S}\left( \widehat{b}\right) \right) }{d\left( b_{\mu },\widehat{b}%
\right) },\text{ for all }\mu \in \mathbf{]}0,1\mathbf{]}.
\end{equation*}

Since $\lim_{\mu \downarrow 0}b_{\mu }=\widehat{b},$ $\lim_{\mu \downarrow
0}x_{\mu }=\widehat{x}$ and$,$ for each $\mu \in \mathbf{]}0,1\mathbf{]},$%
\begin{equation*}
b_{\mu }\in \dom\mathcal{S}_{k}\text{ and }x_{\mu }\in \mathcal{S}_{k}\left(
b_{\mu }\right) =\mathcal{S}\left( b_{\mu }\right) ,
\end{equation*}%
recalling the definition of calmness modulus, we obtain the claimed
inequality 
\begin{equation*}
\frac{d\left( x,\mathcal{S}\left( \widehat{b}\right) \right) }{d\left( b,%
\widehat{b}\right) }=\lim_{\mu \rightarrow 0}\frac{d\left( x_{\mu },\mathcal{%
S}\left( \widehat{b}\right) \right) }{d\left( b_{\mu },\widehat{b}\right) }%
\leq \clm\mathcal{S}\left( \widehat{b},\widehat{x}\right) .
\end{equation*}
\end{dem}

\begin{defn}
\label{Def_subdivision_connecting}Let $\mathcal{S}=\tbigcup\nolimits_{i\in I}%
\mathcal{S}_{i}$ be a wcpc mapping. Let $b,\widehat{b}\in \dom\mathcal{S}$%
.\thinspace\ We say that a subdivision $0=:\mu _{0}<\mu _{1}<...<\mu _{N}:=1$
of $[0,1]$ together with a family of indices $i_{1},...,i_{N}\in I$ connect $%
b$ with $\widehat{b}$ if for all $k\in \{1,...,N\}$ and all $\mu \in [\mu
_{k-1},\mu _{k}]$ there holds 
\begin{equation}
\widehat{b}+\mu \left( b-\widehat{b}\right) \in \dom\mathcal{S}_{i_{k}}
\label{eq_connect_i}
\end{equation}%
(equivalently, $\mathcal{S}\left( \widehat{b}+\mu \left( b-\widehat{b}%
\right) \right) =\mathcal{S}_{i_{k}}\left( \widehat{b}+\mu \left( b-\widehat{%
b}\right) \right) ,$ whenever $\mu \in [\mu _{k-1},\mu _{k}]$).
\end{defn}

\begin{lem}
\label{Lem_subdivision_existence}Let $\mathcal{S}=\tbigcup\nolimits_{i\in I}%
\mathcal{S}_{i}$ be a wcpc mapping. For every pair $b,\widehat{b}\in \dom%
\mathcal{S}$ there are $0=:\mu _{0}<\mu _{1}<...<\mu _{N}:=1,$ with $N\leq
\left\vert I\right\vert $ (the cardinality of $I$) together with a family of
indices $i_{1},...,i_{N}\in I$ connecting $b$ with $\widehat{b}.$
\end{lem}

\begin{dem}
By the convexity of $\limfunc{dom}\mathcal{S},$ the segment $S:=\left\{ 
\widehat{b}+\mu \left( b-\widehat{b}\right) :\mu \in \left[ 0,1\right]
\right\} $ is contained in $\limfunc{dom}\mathcal{S}.$ For each $i\in I,$
let $S_{i}:=S\cap \limfunc{dom}\mathcal{S}_{i}.$ Let $J=\left\{ j\in
I:S_{j}\subsetneqq S_{i}\text{ for some }i\in I\right\} .$ For each $\mu \in %
\left[ 0,1\right] $ let 
\begin{equation*}
ind\left( \mu \right) :=\min \left\{ i\in I\backslash J:\widehat{b}+\mu
\left( b-\widehat{b}\right) \in S_{i}\right\} .
\end{equation*}%
Then, for each $i\in I$, $ind^{-1}\left( i\right) \subset \left[ 0,1\right] $
is either empty of an interval. More in detail, if $\mu _{1},\mu _{2}\in
ind^{-1}\left( i\right) $ with $0\leq \mu _{1}<\mu _{2}\leq 1$ and $%
ind\left( \mu \right) =j<i$ for some $\mu \in \left] \mu _{1},\mu _{2}\right[
,$ then $S_{j}\subsetneqq S_{i}$ because of the convexity of both sets.
Thus, $\left\{ ind^{-1}\left( i\right) :i\in I\backslash J\right\}
\backslash \left\{ \emptyset \right\} $ is a partition of $\left[ 0,1\right] 
$ in at most $\left\vert I\right\vert $ intervals. Let $0=:\mu _{0}<\mu
_{1}<...<\mu _{N}:=1$ be the endpoints of those intervals and for each $%
k=1.,,,.N$ let $i_{k}$ be the unique element of $I\backslash J$ such that $%
\left] \mu _{k-1},\mu _{k}\right[ \subset ind^{-1}\left( i_{k}\right) .$
Observe that the closedness of $S_{i_{k}}$ entails (\ref{eq_connect_i}) for
all $\mu \in [ \mu _{k-1},\mu _{k}].$
\end{dem}

\begin{theo}
\label{Theorem_dist_less_lipsup}Let $\mathcal{S}=\tbigcup\nolimits_{i\in I}%
\mathcal{S}_{i}$ be a wcpc mapping. Let $b,\widehat{b}\in \dom\mathcal{S}$
with $b\neq \widehat{b}$ and consider a subdivision $0=:\mu _{0}<\mu
_{1}<...<\mu _{N}:=1$ together with a family of indices $i_{1},...,i_{N}\in
I $ connecting $b$ with $\widehat{b}.$ Then, for every $x\in \mathcal{S}%
\left( b\right) ,$ there exist points $x^{k}\in \mathcal{S}\left( \widehat{b}%
+\mu _{k}d\right) $ with $k=0,...,N-1$ such that 
\begin{equation*}
\frac{d\left( x,\mathcal{S}\left( \widehat{b}\right) \right) }{d\left( b,%
\widehat{b}\right) }\leq \max \{\clm\mathcal{S}\left( \widehat{b}+\mu
_{k}\left( b-\widehat{b}\right) ,x^{k}\right) \mid k=0,...,N-1\}.
\end{equation*}
\end{theo}

\begin{dem}
We proceed by induction in $N.$ First assume that $N=1.$ Take any $x\in 
\mathcal{S}\left( b\right) $. Then both $\widehat{b}$ and $b$ belong to $\dom%
\mathcal{S}_{i_{1}}.$ Pick $x^{0}\in P_{\mathcal{S}\left( \widehat{b}\right)
}(x)$ and apply Lemma \ref{Lem_ratio_dist} to obtain the aimed inequality%
\begin{equation*}
\frac{d\left( x,\mathcal{S}\left( \widehat{b}\right) \right) }{d\left( b,%
\widehat{b}\right) }\leq \clm\mathcal{S}\left( \widehat{b},x^{0}\right) .
\end{equation*}%
Now, for $N\geq 2,$ consider a subdivision $0=:\mu _{0}<\mu _{1}<...<\mu
_{N}:=1$ together with a family of indices $i_{1},...,i_{N}\in I$ connecting 
$b$ with $\widehat{b}.$ Take any $x\in \mathcal{S}\left( b\right) $ and
consider 
\begin{equation*}
x^{N-1}\in P_{\mathcal{S}_{i_{N}}\left( \widehat{b}+\mu _{N-1}\left( b-%
\widehat{b}\right) \right) }(x)=P_{\mathcal{S}\left( \widehat{b}+\mu
_{N-1}\left( b-\widehat{b}\right) \right) }(x).
\end{equation*}%
On the one hand, Lemma \ref{Lem_ratio_dist} yields

\begin{equation*}
\frac{d\left( x,\mathcal{S}\left( \widehat{b}+\mu _{N-1}\left( b-\widehat{b}%
\right) \right) \right) }{d\left( b,\widehat{b}+\mu _{N-1}\left( b-\widehat{b%
}\right) \right) }\leq \clm\mathcal{S}\left( \widehat{b}+\mu _{N-1}\left( b-%
\widehat{b}\right) ,x^{N-1}\right) ,
\end{equation*}%
equivalently%
\begin{equation}
\left\Vert x-x^{N-1}\right\Vert \leq \clm\mathcal{S}\left( \widehat{b}+\mu
_{N-1}\left( b-\widehat{b}\right) ,x^{N-1}\right) \left\Vert \left( 1-\mu
_{N-1}\right) \left( b-\widehat{b}\right) \right\Vert .  \label{eq_0002}
\end{equation}%
On the other hand, we can write 
\begin{equation}
d\left( x,\mathcal{S}\left( \widehat{b}\right) \right) \leq \left\Vert
x-x^{N-1}\right\Vert +d\left( x^{N-1},\mathcal{S}\left( \widehat{b}\right)
\right) .  \label{eq_0003}
\end{equation}%
Now, consider 
\begin{equation*}
\widetilde{b}:=\widehat{b}+\mu _{N-1}\left( b-\widehat{b}\right)
\end{equation*}%
and subdivision $0=:\widetilde{\mu }_{0}<\widetilde{\mu }_{1}<...<\widetilde{%
\mu }_{N-1}:=1,$ with $\widetilde{\mu }_{i}:=\dfrac{\mu _{i}}{\mu _{N-1}},$ $%
i=1,...,N-1,$ together with the family of indices $i_{1},...,i_{N-1}.$
Observe that this subdivision is connecting $\widehat{b}$ and $\widetilde{b}$
as 
\begin{equation}
\widehat{b}+\widetilde{\mu }_{i}\left( \widetilde{b}-\widehat{b}\right) =%
\widehat{b}+\mu _{i}\left( b-\widehat{b}\right) ,\text{ }i=1,...,N-1.
\label{eq_0003b}
\end{equation}%
Hence, we apply the induction hypothesis with $\widehat{b}$ and $\widetilde{b%
}$ and subdivision $0=:\widetilde{\mu }_{0}<\widetilde{\mu }_{1}<...<%
\widetilde{\mu }_{N-1}:=1$ to conclude that for $x^{N-1}\in \mathcal{S}(%
\widetilde{b})$ we can find points $x^{k}\in \mathcal{S}\left( \widehat{b}+%
\widetilde{\mu }_{k}\left( \widetilde{b}-\widehat{b}\right) \right) $ with $%
k=0,...,N-2$ such that 
\begin{equation}
\frac{d\left( x^{N-1},\mathcal{S}\left( \widehat{b}\right) \right) }{%
\left\Vert \mu _{N-1}\left( b-\widehat{b}\right) \right\Vert }\leq \max \{%
\clm\mathcal{S}\left( \widehat{b}+\widetilde{\mu }_{k}\left( \widetilde{b}-%
\widehat{b}\right) ,x^{k}\right) \mid k=0,...,N-2\}=:\widetilde{\gamma }.
\label{eq_0004}
\end{equation}%
Then, combining (\ref{eq_0002}), (\ref{eq_0003}), (\ref{eq_0003b}) and (\ref%
{eq_0004}) we conclude,%
\begin{eqnarray*}
d\left( x,\mathcal{S}\left( \widehat{b}\right) \right) &\leq &\clm\mathcal{S}%
\left( \widetilde{b},x^{N-1}\right) \left\Vert \left( 1-\mu _{N-1}\right)
\left( b-\widehat{b}\right) \right\Vert +\widetilde{\gamma }\left\Vert \mu
_{N-1}\left( b-\widehat{b}\right) \right\Vert \\
&\leq &\max \{\clm\mathcal{S}\left( \widehat{b}+\mu _{k}\left( b-\widehat{b}%
\right) ,x^{k}\right) \mid k=0,...,N-1\}\left\Vert \left( b-\widehat{b}%
\right) \right\Vert ,
\end{eqnarray*}%
as we wanted to prove.
\end{dem}

\begin{cor}
\label{Cor_Hofmann_Lipusc}Let $\mathcal{S}$ be a $wcpc$ mapping. Then 
\begin{equation*}
\Hof\mathcal{S}=\sup \{\Lipusc\mathcal{S}\left( b\right) \mid b\in \dom%
\mathcal{S}\}=\sup \{\clm\mathcal{S}\left( b,x\right) \mid \left( b,x\right)
\in \gph\mathcal{S}\}.
\end{equation*}
\end{cor}

\begin{dem}
The inequalities \textquotedblleft $\geq $\textquotedblright\ are always
true (recall (\ref{eq_inequalities chain})) and inequalities
\textquotedblleft $\leq $\textquotedblright\ follow from Theorem \ref%
{Theorem_dist_less_lipsup}.
\end{dem}

Now we provide two examples. The first one is intended to illustrate the
previous corollary beyond polyhedral structures, while the second shows that
wcpc assumption of this corollary is not superfluous.

\begin{exa}
\label{Exa1}\emph{Consider} \emph{mapping the set-valued mappings \ }%
\begin{equation*}
\mathcal{S}_{1}(y):=\left\{ 
\begin{array}{l}
\lbrack -1,1],\text{ }\emph{if\ }y\leq 1, \\ 
\emptyset ,\emph{if\ }y>1,%
\end{array}%
\right. \text{ \emph{and} }\mathcal{S}_{2}(y):=\left\{ 
\begin{array}{c}
\emptyset ,\emph{if\ }y<1,, \\ 
\{x\in \mathbb{R\mid }x^{2}\leq y\},\text{\emph{\ if} }y\geq 1.%
\end{array}%
\right. 
\end{equation*}%
\emph{Observe that }$\mathcal{S=S}_{1}\cup \mathcal{S}_{2}$\emph{\ is a wcpc
mapping. We can easily check that} 
\begin{equation*}
\Hof\mathcal{S=}\clm\mathcal{S}\left( 1,1\right) =\left. \frac{d}{dy}\sqrt{y}%
\right\vert _{y=1}=\frac{1}{2}.
\end{equation*}
\end{exa}

\begin{exa}
\label{Exa not unif clm}\emph{Consider} \emph{mappings }$\mathcal{S}_{1},$ $%
\mathcal{S}_{2}:\mathbb{R}\longrightarrow \mathbb{R}$ \emph{given by}%
\begin{equation*}
\mathcal{S}_{1}(y)=\left\{ 
\begin{array}{c}
0\text{ }\emph{if\ }y\leq 0, \\ 
1\text{ }\emph{if\ }y>0,%
\end{array}%
\right. \text{ \emph{and} }\mathcal{S}_{2}(y)=\left\{ 
\begin{array}{c}
0\text{ }\emph{if\ }y<0, \\ 
\lbrack 0,1]\text{ }\emph{if\ }y\geq 0.%
\end{array}%
\right.
\end{equation*}%
\emph{Observe that }%
\begin{equation*}
\clm\mathcal{S}_{1}\left( y,x\right) =0\text{ }\forall \left( y,x\right) \in 
\mathrm{gph}\mathcal{S}_{1},\text{ }\Lipusc\mathcal{S}_{1}\left( 0\right)
=+\infty \text{ },
\end{equation*}%
\emph{and}%
\begin{equation*}
\Lipusc\mathcal{S}_{2}\left( y\right) =0\text{ }\forall y\in \dom\mathcal{S}%
_{2},\text{ \emph{while} }\Hof\mathcal{S}_{2}\mathcal{=+\infty }.
\end{equation*}%
\emph{Clearly} $\mathcal{S}_{1}$ \emph{is not }$\emph{wcpc}$\emph{\ since}%
\textbf{\ }$\mathrm{gph}\mathcal{S}_{1}$ \emph{is not a finite union of
closed convex sets. Regarding }$\mathcal{S}_{2},$ \emph{which can be
represented as} $\mathcal{S}_{2}=\mathcal{S}_{21}\cup \mathcal{S}_{22}$, 
\emph{with} 
\begin{equation*}
\mathrm{gph}\mathcal{S}_{21}=]\mathbb{-}\infty ,0]\times \{0\},\text{ }%
\mathrm{gph}\mathcal{S}_{22}=[0,+\infty \lbrack \times \lbrack 0,1],
\end{equation*}%
\emph{observe that\ condition }$\left( ii\right) $ \emph{in Definition \ref%
{def_well_assembled} fails, since }$0\in \dom\mathcal{S}_{21}\cap \dom%
\mathcal{S}_{22}$ \emph{but }$\mathcal{S}_{21}\left( 0\right) \neq \mathcal{S%
}_{22}\left( 0\right) .$
\end{exa}

The last result of this section refines the second equality of Corollary \ref%
{Cor_Hofmann_Lipusc} and constitutes the version of the second equality of
Theorem \ref{theo: Lip extreme points} adapted to wcpc mappings.

\begin{prop}
Let $\mathcal{S}=\cup _{i\in I}\mathcal{S}_{i}$ be a wcpc mapping. Let $%
\overline{b}\in \dom\mathcal{S}$ and $J:=\left\{ i\in I:\overline{b}\in \dom%
\mathcal{S}_{i}\right\} $. Then 
\begin{equation*}
\Lipusc\mathcal{S}\left( \overline{b}\right) =\max_{j\in J}\sup_{x\in 
\mathcal{S}_{j}\left( \overline{b}\right) }\clm\mathcal{S}_{j}\left( 
\overline{b},x\right)
\end{equation*}
\end{prop}

\begin{dem}
The inequality $\geq $ holds by definition. For the opposite inequality, we
have the following: 
\begin{eqnarray*}
\Lipusc\mathcal{S}\left( \overline{b}\right) &=&\underset{b\rightarrow 
\overline{b}}{\lim \sup }\sup_{x\in \mathcal{S}\left( b\right) }\frac{%
d\left( x,\mathcal{S}\left( \overline{b}\right) \right) }{d\left( b,%
\overline{b}\right) } \\
&\leq &\max_{j\in J}\underset{b\rightarrow \overline{b},~b\in \mathrm{dom}%
\mathcal{S}_{j}}{\lim \sup }\sup_{x\in \mathcal{S}_{j}\left( b\right) }\frac{%
d\left( x,\mathcal{S}_{j}\left( \overline{b}\right) \right) }{d\left( b,%
\overline{b}\right) } \\
&\leq &\max_{j\in J}\underset{b\rightarrow \overline{b},~b\in \mathrm{dom}%
\mathcal{S}_{j}}{\lim \sup }~\sup_{x\in \mathcal{S}_{j}\left( b\right)
}~\sup_{\overline{x}\in P_{\mathcal{S}_{j}\left( \overline{b}\right) }\left(
x\right) }\mathrm{clm}\mathcal{S}_{j}\left( \overline{b},\overline{x}\right)
\\
&\leq &\max_{j\in J}\sup_{x\in \mathcal{S}_{j}\left( \overline{b}\right) }%
\mathrm{clm}\mathcal{S}_{j}\left( \overline{b},x\right) .
\end{eqnarray*}%
The first equality is stated in (\ref{eq_Lipusc}). The second line is
derived from the fact that $\mathcal{S}$ is wcpc. More specifically, we
could write $\Lipusc\mathcal{S}\left( \overline{b}\right) $ as a sequential
limit and, due to the finiteness of $I,$ we could take a subsequence with
all the $b$'s in the same $\mathrm{dom}\mathcal{S}_{j},$ which also contains 
$\overline{b}$ since this set is closed. The third line is deduced from \cite%
[Lemma 2]{CCP21} taking into account that each $\gph\mathcal{S}_{j}$ is
closed and convex. The last inequality is evident.
\end{dem}

\subsection{Hoffman constant under RHS perturbations}

This subsection applies previous results to derive the aimed point-based
formula for $\Hof\mathcal{F}_{\overline{c}}^{op}$ in Theorem \ref{Th upper
bound Hof}. First, we introduce the following lemma, which refers to the
minimal KKT index subsets at $\left( \overline{c},b\right) \in \gph\mathcal{F%
}_{\overline{c}}^{op}$ introduced in Definition \ref{def_KKT_subsets}$.$

\begin{lem}
\label{Prop_Mc} We have that $\mathcal{M}_{\overline{c}}=\bigcup_{b\in \dom%
\mathcal{F}}\mathcal{M}_{\overline{c},b}.$
\end{lem}

\begin{dem}
Inclusion `$\supset $' is evident. On the other hand, for any $D\in \mathcal{%
M}_{\overline{c}}$ we have $\left\{ D\right\} =\mathcal{M}_{\overline{c}%
,b^{D}}=\mathcal{M}_{\overline{c},b^{D}}\left( 0_{n}\right) ,$ where $%
b_{t}^{D}:=$ $0$ if $t\in D$ and $b_{t}^{D}:=$ $1$ if $t\in T\backslash D.$
\end{dem}

\begin{theo}
\label{Th upper bound Hof}Let $-\overline{c}\in \cone\left\{ a_{t},~t\in
T\right\} .\ $One has 
\begin{eqnarray*}
\Hof\mathcal{F}_{\overline{c}}^{op} &=&\max_{\left( b,x\right) \in \gph%
\mathcal{F}_{\overline{c}}^{op}}\clm\mathcal{F}_{\overline{c}}^{op}\left(
b,x\right) \\
&=&\max_{\substack{ D\subset S\subset T  \\ D\in \mathcal{M}_{\overline{c}}}}%
\left\{ d_{\ast }\left( 0_{n},\mathrm{end\,conv}\left\{ a_{t},~t\in
S;~-a_{t},~t\in D\right\} \right) \right\} ^{-1}.
\end{eqnarray*}
\end{theo}

\begin{dem}
From Corollary \ref{Cor_Hofmann_Lipusc}, Theorem \ref{theo: clm feasible}
and Lemma \ref{Prop_Mc}, we have 
\begin{eqnarray}
\Hof\mathcal{F}_{\overline{c}}^{op} &=&\sup_{\left( b,x\right) \in \gph%
\mathcal{F}_{\overline{c}}^{op}}\clm\mathcal{F}_{\overline{c}}^{op}\left(
b,x\right)  \notag \\
&=&\max_{\left( b,x\right) \in \gph\mathcal{F}_{\overline{c}%
}^{op}}\max_{D\in \mathcal{M}_{\overline{c},b}}\left[ d_{\ast }\left( 0_{n},%
\mathrm{end\,conv}\left\{ 
\begin{array}{l}
a_{t},~t\in T_{b}\left( x\right) ; \\ 
-a_{t},~t\in D%
\end{array}%
\right\} \right) \right] ^{-1}  \notag \\
&\leq &\max_{\substack{ D\subset S\subset T  \\ D\in \mathcal{M}_{\overline{c%
}}}}\left\{ d_{\ast }\left( 0_{n},\mathrm{end\,conv}\left\{ a_{t},~t\in
S;~-a_{t},~t\in D\right\} \right) \right\} ^{-1}.  \label{eq_menor}
\end{eqnarray}%
Observe that we have written directly `max' instead of `sup' in the second
expression since the family of minimal KKT subset of indices $\left\{ 
\mathcal{M}_{\overline{c},b}:b\in \dom\mathcal{F}\right\} $ is finite. In
particular, the first equality of the current theorem is already proved.

To prove the converse inequality `$\geq $' in (\ref{eq_menor}) it will be
enough to prove that, for each $D\in \mathcal{M}_{\overline{c}}$ and each $%
D\subset S\subset T,$ there exist $b_{D,S}\in \dom\mathcal{F}$ and $%
x_{D,S}\in \mathcal{F}^{op}\left( \overline{c},b_{D,S}\right) $ such that $%
S=T_{b_{D,S}}\left( x_{D,S}\right) ,$ which automatically implies $D\in 
\mathcal{M}_{\overline{c},b_{D,S}}$ (because $D\in \mathcal{M}_{\overline{c}%
} $). This can be done by just taking 
\begin{equation*}
x_{D,S}=0_{n}\text{ and }b_{D,S}\left( t\right) =\left\{ 
\begin{tabular}{ll}
$0$ & if $t\in S,$ \\ 
$1$ & if $t\in T\backslash S.$%
\end{tabular}%
\right.
\end{equation*}
\end{dem}

\begin{rem}
\emph{An alternative proof for the first equality of the previous theorem
could be given by adapting the proof of \cite[Theorem 2.1]{Li94} to our
multifunction }$\mathcal{F}_{\overline{c}}^{op}$\emph{, whose domain is not} 
$\mathbb{R}^{m}$ \emph{but a polyhedral subset of it. According to the
terminology of that paper (see also Robinson \cite{Rob81}), as a consequence
of its polyhedral character,} $\mathcal{F}_{\overline{c}}^{op}$ \emph{is
locally upper Lipschitz (Lipschitz upper semicontinuous in the sense of (\ref%
{eq_008})) with the same constant for all }$b\in \dom\mathcal{F}_{\overline{c%
}}^{op}$ \emph{(and possible different associated neighborhoods) and
Lipschitz continuous on} $\dom\mathcal{F}_{\overline{c}}^{op}$ \emph{(hence
Hausdorff lower semicontinuous when restricted to its domain)}. \emph{%
Therefore, the adapted version of \cite[Theorem 2.1]{Li94} states that any
locally upper Lipschitz constant of }$\mathcal{F}_{\overline{c}}^{op}$ \emph{%
becomes, in fact, a Lipschitz constant; from this implication one can derive
the equality} 
\begin{equation*}
\Hof\mathcal{F}_{\overline{c}}^{op}=\sup_{b\in \dom\mathcal{F}_{\overline{c}%
}^{op}}\Lipusc\mathcal{F}_{\overline{c}}^{op}\left( b\right) ,
\end{equation*}%
\emph{which coincides with} $\sup_{\left( b,x\right) \in \gph\mathcal{F}_{%
\overline{c}}^{op}}\clm\mathcal{F}_{\overline{c}}^{op}\left( b,x\right) $ 
\emph{as a consequence of Theorem \ref{theo: Lip extreme points}.}

\emph{Observe that the approach of the current paper is slightly more
general as the equality between the Hoffman constant and the supremum of
calmness moduli is established for wcpc mappings}.
\end{rem}

\section{Break steps and directional behavior of the optimal set}

This section is focused on some features of the optimal set mapping $%
\mathcal{F}_{\overline{c}}^{op}$ along the segment determined by two
elements $\overline{b}$ and $b$ of its domain; \ i.e., on the behavior of $%
\mathcal{F}_{\overline{c}}^{op}\left( \overline{b}+\mu \left( b-\overline{b}%
\right) \right) ,$ provided that $\overline{b},b\in \dom\mathcal{F}_{%
\overline{c}}^{op}$ and $\mu \in [0,1]$. Specifically, the main contribution
of this section is to provide a way of computing a subdivision $0=:\mu
_{0}<\mu _{1}<...<\mu _{N}:=1$ together with a family $D_{1},...,D_{N}\in 
\mathcal{M}_{\overline{c}}$ connecting $b$ with $\overline{b}.$ Recall that
the existence of such a subdivision is guaranteed since $\mathcal{F}_{%
\overline{c}}^{op}$ is a $\mathrm{wcpc}$ mapping (recall Lemma \ref%
{Lem_subdivision_existence}).

\begin{lem}
\label{Lem minimals no reapear}Let $\overline{b},b\in \dom\mathcal{F}$.
Assume $\mathcal{M}_{\overline{c},\overline{b}}=\mathcal{M}_{\overline{c}%
,b}\neq \emptyset .$ Then, 
\begin{equation*}
\mathcal{M}_{\overline{c},\overline{b}+\mu \left( b-\overline{b}\right) }=%
\mathcal{M}_{\overline{c},\overline{b}}\text{ for all }\mu \in [0,1].
\end{equation*}
\end{lem}

\begin{dem}
First observe that $\overline{b}+\mu \left( b-\overline{b}\right) \in \dom%
\mathcal{F}$ for all $\mu \in ]0,1[$ because of the convexity of $\dom%
\mathcal{F}.$ Fix arbitrarily $\overline{x}\in \mathcal{F}_{\overline{c}%
}^{op}\left( \overline{b}\right) $ and $x\in \mathcal{F}_{\overline{c}%
}^{op}\left( b\right) $ and consider the convex combination 
\begin{equation*}
\left( b^{\mu },x^{\mu }\right) :=\left( \overline{b},\overline{x}\right)
+\mu \left( b-\overline{b},x-\overline{x}\right) \in \gph\mathcal{F},\text{ }%
\mu \in ]0,1[.
\end{equation*}%
It is immediate from the definitions that 
\begin{equation}
T_{b^{\mu }}\left( x^{\mu }\right) =T_{\overline{b}}\left( \overline{x}%
\right) \cap T_{b}\left( x\right) ,\text{ for all }\mu \in ]0,1[.
\label{eq_active_lambda}
\end{equation}%
More in detail, just observe that, for any $t\in T,$ 
\begin{equation}
a_{t}^{\prime }x^{\mu }-b_{t}^{\mu }=\left( 1-\mu \right) \left(
a_{t}^{\prime }\overline{x}-\overline{b}_{t}\right) +\mu \left(
a_{t}^{\prime }x-b_{t}\right) \leq 0,  \label{eq_axblambda}
\end{equation}%
and equality holds if and only if $a_{t}^{\prime }\overline{x}-\overline{b}%
_{t}=a_{t}^{\prime }x-b_{t}=0.$

Take any $\mu \in \left] 0,1\right[ $ and let us show that $\mathcal{M}_{%
\overline{c},b^{\mu }}=\mathcal{M}_{\overline{c},\overline{b}}$. Observe
that any $D\in \mathcal{M}_{\overline{c},\overline{b}}=\mathcal{M}_{%
\overline{c},b}$ is contained in $T_{b^{\mu }}\left( x^{\mu }\right) $
because of (\ref{eq_active_lambda}), and clearly $D$ is minimal among the
subsets of $T_{b^{\mu }}\left( x^{\mu }\right) $ satisfying $-\overline{c}%
\in \cone\{a_{t},\,t\in D\}$, since this minimality happens, for instance,
in the subsets of $T_{\overline{b}}\left( \overline{x}\right) .$ In other
words, $D\in \mathcal{M}_{\overline{c},b^{\mu }}$.

In order to check the converse inclusion, take any $D\in \mathcal{M}_{%
\overline{c},b^{\mu }},$ in particular $D\in T_{b^{\mu }}\left( x^{\mu
}\right) $ and hence $D\in T_{\overline{b}}\left( \overline{x}\right) .$
Moreover, the minimality of $D$ among the subsets of $T_{\overline{b}}\left( 
\overline{x}\right) $ is derived from the minimality over the subsets of $%
T_{b^{\mu }}\left( x^{\mu }\right) .$
\end{dem}

\begin{defn}
Given $\overline{b},b\in \dom\mathcal{F}$, $\overline{b}\neq b,$ we define
the \emph{break step set between }$\overline{b}$ \emph{and }$b$\emph{\ }by%
\emph{:} 
\begin{equation*}
\mathcal{B}\left( \overline{b},b\right) :=\left\{ \mu \in ]0,1[\left\vert
\exists \nu _{r}\rightarrow \mu \text{ , }\mathcal{M}_{c,\overline{b}+\nu
_{r}\left( b-\overline{b}\right) }\subsetneqq \mathcal{M}_{c,\overline{b}%
+\mu \left( b-\overline{b}\right) }\text{ }\forall r\in \mathbb{N}\right.
\right\} .
\end{equation*}
\end{defn}

\begin{rem}
\label{Rem_00}\emph{Because of Lemma \ref{lem: tech}, in the previous
definitions we might replace }"$\subsetneqq $"\emph{\ with }"$\neq $".
\end{rem}

\begin{prop}
\label{Prop_Bpp}Given $\overline{b},b\in \dom\mathcal{F}$, $\overline{b}\neq
b,$ we have:

$\left( i\right) $ $\mathcal{B}\left( \overline{b},b\right) $ is a finite
set (possibly empty).

$\left( ii\right) $ If $\mathcal{B}\left( \overline{b},b\right) =\{\mu
_{1},\mu _{2},...,\mu _{N}\},$ with $0=:\mu _{0}<\mu _{1}<\mu _{2}<...<\mu
_{N}<\mu _{N+1}:=1,$ then 
\begin{equation*}
\mu \mapsto \mathcal{M}_{c,\overline{b}+\mu \left( b-\overline{b}\right) }%
\text{ is constant on }]\mu _{i-1},\mu _{i}[\text{ for any }i=1,2,...,N+1.
\end{equation*}%
(Take $N=0$ when $\mathcal{B}\left( \overline{b},b\right) =\emptyset ,$ in
which case we just have $\mathcal{M}_{\overline{c},\overline{b}+\mu \left( b-%
\overline{b}\right) }=\mathcal{M}_{\overline{c},\overline{b}}$ for all $\mu
\in [0,1]$).

$\left( iii\right) $ Let $\mathcal{B}\left( \overline{b},b\right) =\{\mu
_{1},\mu _{2},...,\mu _{N}\},$ with $0:=\mu _{0}<\mu _{1}<...<\mu _{N}<\mu
_{N+1}:=1$ and, for each $i=1,...,N+1,$ fix any 
\begin{equation*}
D_{i}\in \mathcal{M}_{\overline{c},\overline{b}+\mu \left( b-\overline{b}%
\right) },\text{ }\mu \in ]\mu _{i-1},\mu _{i}[.
\end{equation*}%
Then the subdivision $0:=\mu _{0}<\mu _{1}...\mu _{N}<\mu _{N+1}=1$ together
with the family $D_{1},...,D_{N+1}\in \mathcal{M}_{\overline{c}}$ connect $b$
with $\overline{b}.$
\end{prop}

\begin{dem}
$\left( i\right) $ Reasoning by contradiction, assume that $\mathcal{B}%
\left( \overline{b},b\right) $ is infinite, take a sequence of scalars $%
\{\mu _{r}\}_{r\in \mathbb{N}}\subset \mathcal{B}\left( \overline{b}%
,b\right) $ and consider the corresponding sequence of subsets $\left\{ 
\mathcal{M}_{\overline{c},\overline{b}+\mu _{r}\left( b-\overline{b}\right)
}\right\} _{r\in \mathbb{N}}.$ Since $\mathcal{M}_{\overline{c},\overline{b}%
+\mu _{r}\left( b-\overline{b}\right) }\subset \mathcal{M}_{\overline{c}},$
for all $r,$ and $\mathcal{M}_{\overline{c}}$ is finite, there exist
infinitely many repeated subsets in $\left\{ \mathcal{M}_{\overline{c},%
\overline{b}+\mu _{r}\left( b-\overline{b}\right) }\right\} _{r\in \mathbb{N}%
}.$ In particular, we can take three break steps $\mu _{r_{1}}<\mu
_{r_{2}}<\mu _{r_{3}}$ with the same minimal KKT set of indices. Then, since 
$\mathcal{M}_{\overline{c},\overline{b}+\mu _{r_{1}}\left( b-\overline{b}%
\right) }=\mathcal{M}_{\overline{c},\overline{b}+\mu _{r_{3}}\left( b-%
\overline{b}\right) },$ applying the previous lemma (with $\overline{b}+\mu
_{r_{1}}\left( b-\overline{b}\right) $ and $\overline{b}+\mu _{r_{3}}\left(
b-\overline{b}\right) $ playing the roles of $\overline{b}$ and $b,$
respectively), $\mathcal{M}_{\overline{c},\overline{b}+\mu \left( b-%
\overline{b}\right) }=$\linebreak $\mathcal{M}_{\overline{c},\overline{b}%
+\mu _{r_{1}}\left( b-\overline{b}\right) }$ for all $\mu \in [ \mu
_{r_{1}},\mu _{r_{3}}],$ which, according to Lemma \ref{Lem minimals no
reapear}, contradicts the fact that $\mu _{r_{2}}$ is also a break step
between $\overline{b}$ and $b.$

$\left( ii\right) $ Fix any $i\in \{1,...,N+1\}$ and let us see that $%
\mathcal{M}_{c,\overline{b}+\mu \left( b-\overline{b}\right) }$ is constant
on $]\mu _{i-1},\mu _{i}[.$ Arguing by contradiction, assume that there
exist $\overline{\mu }$ and $\widetilde{\mu }$ with $\mu _{i-1}<\overline{%
\mu }<\widetilde{\mu }<\mu _{i}$ such that $\mathcal{M}_{\overline{c},%
\overline{b}+\overline{\mu }\left( b-\overline{b}\right) }\neq \mathcal{M}_{%
\overline{c},\overline{b}+\widetilde{\mu }\left( b-\overline{b}\right) }.$
Define 
\begin{equation*}
\alpha :=\sup \left\{ \mu >\mu _{i-1}\mid \mathcal{M}_{\overline{c},%
\overline{b}+\mu \left( b-\overline{b}\right) }=\mathcal{M}_{\overline{c},%
\overline{b}+\overline{\mu }\left( b-\overline{b}\right) }\right\} .
\end{equation*}%
Observe that $\alpha >\overline{\mu }$ since $\overline{\mu }\notin \mathcal{%
B}\left( \overline{b},b\right) $ (which is deduced from the definition of
break step together with Remark \ref{Rem_00}). Indeed, $\alpha \geq \mu _{i}$
(contradicting the existence of $\widetilde{\mu }$) since, if we had $\alpha
<\mu _{i}$ we would attain a contradiction by distinguishing two cases:

\emph{Case 1}$:$ If $\mathcal{M}_{\overline{c},\overline{b}+\alpha \left( b-%
\overline{b}\right) }=\mathcal{M}_{\overline{c},\overline{b}+\overline{\mu }%
\left( b-\overline{b}\right) }$ by definition of supremum, there would exist
a decreasing sequence of scalars $\nu _{j}\downarrow \alpha $ with $\mathcal{%
M}_{\overline{c},\overline{b}+\nu _{j}\left( b-\overline{b}\right) }\neq 
\mathcal{M}_{\overline{c},\overline{b}+\alpha \left( b-\overline{b}\right)
}. $ Then we would have $\alpha \in \mathcal{B}\left( \overline{b},b\right) $
which represents a contradiction (observe that $\mu _{i-1}<\alpha <\mu _{i}$%
).

\emph{Case 2}$:$ If $\mathcal{M}_{\overline{c},\overline{b}+\alpha \left( b-%
\overline{b}\right) }\neq \mathcal{M}_{\overline{c},\overline{b}+\overline{%
\mu }\left( b-\overline{b}\right) },$ again by the definition of supremum,
there would exists an increasing sequence $\nu _{j}\uparrow \alpha $ with $%
\mathcal{M}_{\overline{c},\overline{b}+\nu _{j}\left( b-\overline{b}\right)
}\neq \mathcal{M}_{\overline{c},\overline{b}+\alpha \left( b-\overline{b}%
\right) }.$ Again, we would attain the contradiction $\alpha \in \mathcal{B}%
\left( \overline{b},b\right) $.

$\left( iii\right) $ is a direct consequence of $\left( i\right) $ and $%
\left( ii\right) $ together with Proposition \ref{Prop_FOp_well_connected}.
Specifically, $\mathcal{F}_{\overline{c}}^{op}=\bigcup_{D\in \mathcal{M}_{%
\overline{c}}}\mathcal{P}_{D},$ subsets $D_{1},...,D_{N}\in \mathcal{M}_{%
\overline{c}}$ and 
\begin{equation*}
\overline{b}+\mu \left( b-\overline{b}\right) \in \dom\mathcal{P}_{D_{i}}%
\text{ for all }\mu \in \lbrack \mu _{i-1},\mu _{i}],
\end{equation*}%
since $D_{i}\in \mathcal{M}_{\overline{c},\overline{b}+\mu \left( b-%
\overline{b}\right) }$ for all $\mu \in ]\mu _{i-1},\mu _{i}[,$ and taking
the closedness of $\dom\mathcal{P}_{D_{i}}$ into account (recall Definition %
\ref{Def_subdivision_connecting}).
\end{dem}

The following result is a direct consequence of Theorem \ref%
{Theorem_dist_less_lipsup} and Propositions \ref{Prop_FOp_well_connected}
and \ref{Prop_Bpp}.

\begin{prop}
\label{Theorem_sect4}Let $b,\overline{b}\in \dom\mathcal{F},$ $b\neq 
\overline{b},$ and consider the set of break steps $\mathcal{B}\left( 
\overline{b},b\right) =\{\mu _{1},\mu _{2},...,\mu _{N}\},$ with $0:=\mu
_{0}<\mu _{1}<\mu _{2}<...<\mu _{N}<\mu _{N+1}:=1.$ Then for every $x\in 
\mathcal{F}_{\overline{c}}^{op}\left( b\right) $ one has 
\begin{equation*}
\frac{d\left( x,\mathcal{F}_{\overline{c}}^{op}\left( \overline{b}\right)
\right) }{d_{\infty }\left( b,\overline{b}\right) }\leq \max \{\Lipusc%
\mathcal{F}_{\overline{c}}^{op}\left( \overline{b}+\mu _{k}d\right) \mid
k=0,...,N\}.
\end{equation*}
\end{prop}

\section{Conclusions, perspectives and examples}

One of the main contributions of the current work is to ensure the
fulfilment of the following equalities 
\begin{equation}
\Hof\mathcal{S}=\sup \{\Lipusc\mathcal{S}\left( b\right) \mid b\in \dom%
\mathcal{S}\}=\sup \{\clm\mathcal{S}\left( b,x\right) \mid \left( b,x\right)
\in \gph\mathcal{S}\},  \label{eq_01}
\end{equation}%
provided that $\mathcal{S}$ is a well-connected piecewise convex mapping,
which is the case of the argmin mapping $\mathcal{F}_{\overline{c}}^{op}$
(under RHS perturbations). The global Lipschitzian behavior of optimal
solutions is characterized through the local one. On the other hand, in the
particular case when $\mathcal{S}=\mathcal{F}_{\overline{c}}^{op}$ we can go
further to derive a point-based expression for such constant:%
\begin{equation}
\Hof\mathcal{F}_{\overline{c}}^{op}=\max_{\substack{ D\subset S\subset T  \\ %
D\in \mathcal{M}_{\overline{c}}}}\left\{ d_{\ast }\left( 0_{n},\mathrm{%
end\,conv}\left\{ a_{t},~t\in S;~-a_{t},~t\in D\right\} \right) \right\}
^{-1},  \label{eq_00002}
\end{equation}%
where $\mathcal{M}_{\overline{c}}$ is defined in (\ref{eq_002}).

The remaining part of this section is divided into two subsections devoted
to related antecedents and future research.

\subsection{Comparison with background}

We can find a variety of different Lipschitz constants (upper estimates for
the Hoffman constant) of the argmin mapping of linear programs in the
literature (cf. \cite[Section 6]{HoTu02}, \cite{Li93}). In \cite[Section 6]%
{HoTu02} such a bound is given for the argmin mapping of the dual program to
(\ref{eq_LPproblem}), which depends only on the matrix $A$ and, therefore,
it may not be exact. In contrast, formula (\ref{eq_00002}) provides the
exact Hoffman constant, which depends not only on $A$, but also on the
coefficient vector of the objective function.

The following example shows that Lipschitz constant of $\mathcal{F}_{%
\overline{c}}^{op}$ given in \cite{Li93}, which does involve $\overline{c},$
can also be strictly greater than $\Hof\mathcal{F}_{\overline{c}}^{op}..$

\begin{exa}
\label{Exa_Li}\emph{Consider the argmin mapping }$\mathcal{F}_{\overline{c}%
}^{op}:\mathbb{R}^{4}\rightrightarrows \mathbb{R}^{2},$ \emph{with }$%
\overline{c}=\left( 1,0\right) ^{\prime },$ \emph{given by} 
\begin{equation*}
\mathcal{F}_{\overline{c}}^{op}\left( b\right) :=\arg \min \{x_{1}\mid
x_{1}\leq b_{1},x_{2}\leq b_{2},-2x_{1}\leq b_{3},-x_{2}\leq b_{4}\}
\end{equation*}%
\emph{and assume that} $\mathbb{R}^{2}$ \emph{is endowed with the Euclidean
norm. Here we have }$\mathcal{M}_{\overline{c}}=\{\{3\}\}$ \emph{and,
applying Theorem \ref{Th upper bound Hof}, one easily checks that the
Hoffman constant can be attained at} $S=\{2,3\}$ \emph{yielding} 
\begin{equation*}
\Hof\mathcal{F}_{\overline{c}}^{op}=\left\{ d_{\ast }\left( 0_{2},\mathrm{%
end\,conv}\left\{ a_{2,},a_{3},-a_{2}\right\} \right) \right\} ^{-1}=\sqrt{5}%
/2.
\end{equation*}%
\emph{Concerning Li's constant introduced in \cite[Formula (1.5) ]{Li93}
(stated as Lipschitz constant in Theorem 2.5 of the same paper), in this
example it writes (in the terminology used in \cite{Li93}) as}%
\begin{equation*}
\gamma _{\infty ,2}\left( A,\overline{c}\right) =\max_{I\in W\left( A\right)
}\sup \left\{ d_{2}\left( A_{I,0}^{+}u,\left\{ x\in \mathbb{R}^{2}\mid
Ax\leq 0_{4},\overline{c}^{\prime }x=0\right\} \right) :\left\Vert
u\right\Vert _{\infty }=1\right\} ,
\end{equation*}%
\emph{where} $A$\emph{\ is the matrix formed by the four rows} $%
a_{t}^{\prime },$ $t=1,...4,$ $W\left( A\right) =\{I\mid \mathrm{rank}A=%
\mathrm{rank}A_{I}=\left\vert I\right\vert \},$ $A_{I}$ \emph{is the matrix
formed by the rows }$a_{t}^{\prime },$ $t\in I,$ $A_{I,0}\,$\emph{is
obtained by replacing rows} $a_{t}^{\prime }$ \emph{of} $A$\emph{\ for} $%
t\notin I$ \emph{by} $0_{2}^{\prime }$ \emph{and} $A_{I,0}^{+}$ \emph{is the
Moore-Penrose inverse of }$A_{I,0}.$ \emph{Obviously} $\left\{ x\in \mathbb{R%
}^{2}\mid Ax\leq 0_{4},\overline{c}^{\prime }x=0\right\} =\{0_{2}\},$ $%
W\left( A\right) =\{\{1,2\},\{1,4\},\{2,3\},\{3,4\}\}$ \emph{and one easily
checks}%
\begin{equation*}
\gamma _{\infty ,2}\left( A,\overline{c}\right) =\max_{I\in W\left( A\right)
}\sup \left\{ \left\Vert A_{I,0}^{+}u\right\Vert :\left\Vert u\right\Vert
_{\infty }=1\right\} =\sqrt{2}>\Hof\mathcal{F}_{\overline{c}}^{op}.
\end{equation*}
\end{exa}

Hoffman constant of the feasible set mapping has been more intensively
studied than that of the argmin mapping. For the sake of completeness let us
include some paragraphs about condition measures and their connection with $%
\Hof\mathcal{F}$. Here we do not use dimensions $m$ and $n$ in order to
avoid possible misunderstandings. For a matrix $M\in \mathbb{R}^{k\times l}$
and norms $\left\Vert \cdot \right\Vert _{\alpha }$ in $\mathbb{R}^{l}$ and $%
\left\Vert \cdot \right\Vert _{\beta }$ in $\mathbb{R}^{k},$ we consider the
induced matrix norm 
\begin{equation*}
\left\Vert M\right\Vert _{\alpha \beta }=\max \left\{ \left\Vert
Mx\right\Vert _{\beta }\mid \left\Vert x\right\Vert _{\alpha }\leq 1\right\}
.
\end{equation*}%
Any matrix $M$, and more specifically its kernel, $\ker M,$ has an
associated \emph{linear matroid}, as described in the survey paper \cite[%
Section 1]{ENV22}. Matroid theory combines in an abstract axiomatic
framework some of the essential features of linear algebra and graph theory,
and borrows some terminology from both areas, as independence, bases, and
circuits, among others. The reader is addressed to \cite{ENV22} for more
specific details. The matroid setting constitutes a powerful context for
dealing with condition measures of matrices and diameters of polyhedra, and
these measures constitute important ingredients in the complexity analysis
of certain algorithms, as interior-point methods for solving linear programs
in polynomial time; see \cite{HoTu02}. We focus our attention on the
condition measures $\bar{\chi}\left( M\right) $ and $\chi \left( M\right) $
considered in \cite[Section 2]{HoTu02}. The first one is closely related to
the \emph{fractional circuit imbalance measure }of $\ker M$ given in \cite[%
Definition 1.1]{ENV22}, denoted by $\kappa _{M}$. There is a close
relationship between these measures and the Hoffman constant of the feasible
set mapping $\mathcal{F}_{M}$ defined as in (\ref{eq_F(b)}) with $%
A=M^{\prime }$ (the transpose of $M$). This relationship is quantified in 
\cite[Section 4]{HoTu02}, paying special attention to the case when both $%
\mathbb{R}^{k}$ and $\mathbb{R}^{l}$ are endowed with the respective
Euclidean norms. Let us mention that \cite[Section 5]{HoTu02} deals with the
relationship between $\Hof\mathcal{F}_{M^{\prime }}$ and Ye's complexity
measure for primal and dual linear programs.

Assume that $M\in \mathbb{R}^{k\times l}$ is a full row rank matrix, denote
by $M_{B}$ the submatrix of $M$ whose \emph{columns }have indices in $%
B\subset \left\{ 1,...,l\right\} $ and let $\mathcal{B}\left( M\right) $ the
family of \emph{bases} of $M$ (recall that $B\in \mathcal{B}\left( M\right)
:\Leftrightarrow B$ is formed by $k$ elements and $M_{B}$ is a square
nonsingular matrix). Both $\bar{\chi}\left( M\right) $ and $\kappa _{M}$ are
given by the same expression%
\begin{equation*}
\max \left\{ \left\Vert M_{B}^{-1}M\right\Vert _{\alpha \beta }\mid B\in 
\mathcal{B}\left( M\right) \right\} ,
\end{equation*}%
with the only difference of the choice of norms. For $\bar{\chi}\left(
M\right) $ both norms, $\left\Vert \cdot \right\Vert _{\alpha }$ in $\mathbb{%
R}^{l}$ and $\left\Vert \cdot \right\Vert _{\beta }$ in $\mathbb{R}^{k},$
are the Euclidean norms, see \cite[Proposition 2.3]{HoTu02} or \cite[%
Proposition 4.1]{ENV22}, whereas for $\kappa _{M}$ both are the maximum
norms $\left\Vert \cdot \right\Vert _{\infty }$, \cite[Proposition 3.1]%
{ENV22}. The following expression for $\chi \left( M\right) ,$ using the
same ingredients, can be found in \cite[Formula (9)]{HoTu02}:%
\begin{equation*}
\chi \left( M\right) =\max \left\{ \left\Vert M_{B}^{-1}\right\Vert _{2}\mid
B\in \mathcal{B}\left( M\right) \right\} ,
\end{equation*}%
where the subscript 2 means that the Euclidean norm is considered in both $%
\mathbb{R}^{l}$ and $\mathbb{R}^{k}.$ In order to compare the previous
expression with the Hoffman constant of $\mathcal{F}_{M^{\prime }},$ as a
consequence of \cite[Theorem 5]{CCP21}, we can write 
\begin{equation*}
\Hof\mathcal{F}_{M^{\prime }}=\max \left\{ \left\Vert \left( M_{B}^{\prime
}\right) ^{-1}\right\Vert _{\alpha \beta }\mid B\in \mathcal{B}\left(
M\right) \right\} ,
\end{equation*}%
$\left\Vert \cdot \right\Vert _{\alpha }$ being the maximum norm and $%
\left\Vert \cdot \right\Vert _{\beta }$ an arbitrary norm.

\subsection{Future research}

There is another type of Hoffman constants located at a given (nominal)
point of the domain of the multifunction under consideration. Following the
terminology of \cite{CCP21}, the \emph{Hoffman modulus} of $\mathcal{S}$ 
\emph{at} $\overline{b}\in \dom\mathcal{S}$ is defined as 
\begin{equation}
\Hof\mathcal{S}(\overline{b}):=\sup_{\left( b,x\right) \in \gph\mathcal{S}}%
\dfrac{d(x,\mathcal{S}(\overline{b}))}{d\left( b,\overline{b}\right) }.
\label{eq_02}
\end{equation}%
It is clear that we can add the following equality to (\ref{eq_01}):%
\begin{equation*}
\Hof\mathcal{S=}\sup_{b\in \dom\mathcal{S}}\Hof\mathcal{S}(b).
\end{equation*}%
In the case when $\gph\mathcal{S}$ is a convex polyhedral set $\Hof\mathcal{S%
}(\overline{b})=\Lipusc\mathcal{S}\left( \overline{b}\right) $ (see \cite[%
Theorem 4 ]{CCP21} for a more general framework where this equality holds)$,$
which is the case of the feasible set mapping $\mathcal{S=F}\emph{.}$
However, this is not the case of a general well-connected piecewise convex
(even polyhedral) mapping, as the following simple example shows.

\begin{exa}
\emph{Consider the} \emph{mapping }$\mathcal{S}:\mathbb{R}\longrightarrow 
\mathbb{R}$ \emph{given by}%
\begin{equation*}
\mathcal{S}(y)=\left\{ 
\begin{array}{c}
0\text{ }\emph{if\ }y\leq 0 \\ 
y\text{ }\emph{if\ }y>0%
\end{array}%
\right. \text{ }
\end{equation*}%
\emph{Observe that }%
\begin{equation*}
\Lipusc\mathcal{S}\left( -1\right) =0\text{ while }\Hof\mathcal{S}(-1)=\Hof%
\mathcal{S=}1.
\end{equation*}

\emph{Another feature about }$\Hof\mathcal{S}(\overline{b})$\emph{\ is that
the supremum in (\ref{eq_02}) may be attained or not. In the previous
example, it is not attained, while, for instance, for the new mapping given
by }$\widetilde{\mathcal{S}}(y)=0,$ $\emph{if\ }y\leq 0,$ $\widetilde{%
\mathcal{S}}(y)=y,$ \emph{if }$y\in [0,1],$ $\widetilde{\mathcal{S}}(y)=1,$%
\emph{\ if }$y\geq 1,$\emph{\ we have} 
\begin{equation*}
\Hof\widetilde{\mathcal{S}}=1>\Hof\widetilde{\mathcal{S}}(-1)=\frac{%
\left\vert \mathcal{S}\left( 1\right) -\mathcal{S}\left( -1\right)
\right\vert }{1-(-1)}=\frac{1}{2}>\Lipusc\widetilde{\mathcal{S}}\left(
-1\right) =0.
\end{equation*}%
\emph{Moreover, in contrast to Lipschitz upper semicontinuity moduli or to
the Hoffman constant, the Hoffman modulus }$\Hof\widetilde{\mathcal{S}}(-1)$%
\emph{\ cannot be expressed in terms of calmness moduli, since the only
possible calmness moduli in this example are either }$0$\emph{\ or }$1$\emph{%
.}
\end{exa}

All the previous situations may also happen in general for mapping $\mathcal{%
F}_{\overline{c}}^{op}.$ Specifically, $\Hof\mathcal{F}_{\overline{c}}^{op}(%
\overline{b})$ may be strictly in between $\Lipusc\mathcal{F}_{\overline{c}%
}^{op}\left( \overline{b}\right) $ and $\Hof\mathcal{F}_{\overline{c}}^{op}$%
, and the supremum defining $\Hof\mathcal{F}_{\overline{c}}^{op}(\overline{b}%
)$ may be attained or not. Moreover, calmness moduli are not enough to
express $\Hof\mathcal{F}_{\overline{c}}^{op}(\overline{b}).$

While we already have (Theorems \ref{theo: clm feasible} and \ref{theo: Lip
extreme points}) point-based formulae (only involving the nominal data) for
computing $\Lipusc\mathcal{F}_{\overline{c}}^{op}\left( \overline{b}\right) $
and $\clm\mathcal{F}_{\overline{c}}^{op}\left( \overline{b},x\right) ,$ for $%
x\in \mathcal{F}_{\overline{c}}^{op}\left( \overline{b}\right) ,$ the
computation of $\Hof\mathcal{F}_{\overline{c}}^{op}(\overline{b})$ remains
as open problem.

\textbf{Acknowledgement: }The authors are indebted to Prof. Diethard Klatte
for his very careful revision of the original version, constructive critical
comments and invaluable help to better integrate the new contributions into
the existing literature. We are also very thankful to the anonymous
referees, whose comments and criticisms have helped us to improve and
reorganize the original version, and also to include additional material and
background references.

\end{document}